\documentclass[11pt]{article}
\usepackage{graphicx}
\usepackage{pgfplots}
\usepackage[dvipsnames]{xcolor}
\usepackage[colorlinks=true,citecolor=blue]{hyperref}
\usepackage{natbib}
\usepackage{graphicx}
\usepackage{amsfonts}
\usepackage{amsmath}
\usepackage{amssymb}
\usepackage{url}
\usepackage{fancyhdr}
\usepackage{indentfirst}
\usepackage{enumerate}
\usepackage{titlesec}
\usepackage{amsthm}
\usepackage{dsfont}
\usepackage[misc]{ifsym}
\usepackage{enumitem}

\theoremstyle{definition}

\newtheorem{example}{Example}

\theoremstyle{plain}
\newtheorem{theorem}{Theorem}
\newtheorem{lemma}{Lemma}
\newtheorem{proposition}{Proposition}
\newtheorem{corollary}{Corollary}
\theoremstyle{remark}
\newtheorem{remark}{Remark}

\usepackage{booktabs}

\usepackage{cases}

\theoremstyle{definition}

\def\N{\mathbb{N}}

\def\p{\mathbb{P}}

\def\R{\mathbb{R}}

\def\B{\mathcal{B}}

\def\ukappa{\overline{\kappa}}
\def\lkappa{\underline{\kappa}}
\def\uC{\overline{\mathcal{C}}}
\def\lC{\underline{\mathcal{C}}}
\def\uJ{\overline{\mathcal{J}}}
\def\lJ{\underline{\mathcal{J}}}
\def\uK{\overline{\mathcal{K}}}
\def\lK{\underline{\mathcal{K}}}
\def\X{\mathcal{X}}

\def\X{\mathcal{X}}

\def\d{\mathrm{d}}

\newcommand{\VaR}{\mathrm{VaR}}

\DeclareMathOperator*{\argmax}{arg\,max}

\usepackage[onehalfspacing]{setspace}

\usepackage{bm}
\usepackage{tikz-qtree}
\usepackage{tikz}

\def\id{\mathds{1}}

\pgfplotsset{compat=1.18}

\title{The upper copula capacity}

\author{Benjamin C\^ot\'e\thanks{Department of Statistics and Actuarial Science, University of Waterloo, Canada.   \texttt{b3cote@uwaterloo.ca}}\and  Ruodu Wang\thanks{Department of Statistics and Actuarial Science, University of Waterloo, Canada.   \texttt{wang@uwaterloo.ca}}}

\date{\today}

\begin{document}
	\maketitle
	\begin{abstract}

Many problems in statistics, combinatorics, econometrics, and other fields are expressible in terms of the upper probability of all copulas.  
We formally study this capacity, which we call the \textit{upper copula} capacity, as a mathematical object. In particular, we examine the capacity of given sets and integration with respect to this capacity. For two-dimensional convex sets, the upper copula capacity equals the smaller of the two coordinate projection lengths. We also give an exact formula for the capacity of lower-tail events for positively weighted sums in any dimension and characterize when the upper Choquet integral equals the supremum of expectations over all copulas. Various problems in different fields can be formulated via the upper-copula capacity, and we discuss how known and new results connect to offer new perspectives on these problems. 

\medskip

\noindent
\textbf{Keywords}: Upper probability, dependence uncertainty,  multipartite fractional matching, ambiguity, joint mixability. 
	\end{abstract}

\section{Introduction}
  

 For $n\in\mathbb{N}$ and a Borel set $A\subseteq [0,1]^n$, let 
\begin{equation*}
    \uC_n(A) = \sup_{C\in\mathbb{C}_n}C(A)\quad \text{and} \quad \lC_n(A)= \inf_{C\in\mathbb{C}_n} C(A),
\end{equation*}
where $\mathbb{C}_n$ is the set of copula measures, i.e., measures with uniform marginals, on $[0,1]^n$.
A supremum (resp.~infimum) of probability measures is a capacity, and it is called an upper (resp.~lower) probability; see e.g.,~\cite{HR09}. Correspondingly, we will refer to $\uC_n$ as the \textit{upper copula capacity} and to $\lC_n$ as the \textit{lower copula capacity}. 
By Sklar's theorem \citep{S59}, any distribution can be decomposed into two parts: its marginals and a copula, expressing its dependence scheme. The upper and lower copula capacities thus describe a problem of dependence uncertainty: the supremal and infimal probability mass attributable to a set under any dependence scheme. 
Upper probabilities, and capacities in general, have been standard objects to represent model ambiguity, for instance in robust statistics \citep{D08, HS73, WF82,HR09}, social sciences \citep{J91} and econometrics \citep{MM04}.

The upper copula capacity differs from the usual upper probabilities considered in this strand of research by the fact that it is not submodular (also called 2-alternating). The set of admissible distributions used in defining an upper probability needs to satisfy some conditions to yield a submodular capacity \citep{HS73,B82}, and the construction method for these sets in the seminal paper by \cite{D08} ensures the resulting upper probability to be 2-alternating. Our set of admissible distributions is rather rooted in the dependence uncertainty framework and this leads to  different construction and results. Similarly, the lower copula capacity is not supermodular, as would otherwise be usual for a lower probability. 
As we will discuss, the upper and lower copula capacities encompass many significant problems related to the existence of dependence scheme, such as joint mixability \citep{WW16} or the merging of $p$-values \citep{VWW22}, and even some problems for which the relation to dependence uncertainty is not directly apparent, 
plank-covering problems \citep{B51,AHKM02} and resource allocation games   \citep{R06}, for instance. The upper copula capacity also relates to the problem of existence of distributions with given marginals \citep{S65}, and we derived one of our results from \cite{K84} who studied this problem.

To reconcile all these problems under the evaluation of the upper or lower copula capacity allows to connect results from different parts of the literature, and more importantly, gives a broader view of the state of the literature on these questions. 

The paper is organized as follows. The mathematical setting of the paper is presented in Section~\ref{sect:preliminaries}. In Section~\ref{sect:uppercopula}, we formally present the upper copula capacity and examine some of its properties. In Section~\ref{sect:geometry}, we discuss the measure of sets with respect to the upper copula capacity, and specifically how one can rely on some geometric characteristics, such as convexity or decreasingness, to compute its value. In Section~\ref{sect:INT}, we discuss Choquet integration with respect to this capacity and how this relates to the Choquet expected utility framework of \cite{S89}. In Section~\ref{sect:app}, we present some of the many problems expressible in terms of the upper copula capacity. Section~\ref{sect:conclusion-copula} concludes.

\section{Capacities and copulas}
\label{sect:preliminaries}

Let $[n]=\{1,\dots,n\}$ and, unless stated otherwise, assume $n\in\mathbb{N}$.   For an $n$-dimensional vector $\mathbf{u}$, we write $u_i$ to denote its $i$th component.  
%
%
%
%
%
Consider the measurable space $([0,1]^{n}, \mathcal{B}_n)$ where $\mathcal{B}_n$ is the Borel $\sigma$-algebra on $[0,1]^{n}$. We write $\mathcal{B}=\mathcal{B}_1$. 
Let $p_i$ be the usual projection of a set $A\in\mathcal{B}_n$ to its $i$th dimension, i.e.~$p_i(A) = \{u_i : \mathbf{u}\in A\}$, for $i\in[n]$. Let $p_i^{-1}$ mean $p_i^{-1}(A) = [0,1]^{i-1}\times A \times [0,1]^{n-i}$ for $A\in\mathcal{B}$, $i\in[n]$. 
For a distribution function $F$, let $F^{-1}$ be its left quantile function.
Throughout, we write $\ell_n$ for the Lebesgue measure on $[0,1]^n$; let $\ell = \ell_1$. We employ the usual topology on $\mathbb{R}^n$ when referring to the openness and the closedness of a set, or when taking its interior or its closure.

  Let $\mathcal F$ be a $\sigma$-field on a set $\Omega$. 
  A \emph{capacity} is a set function $\nu: \mathcal F\to \R_+$ that is \emph{increasing}, meaning $\nu(A)\le \nu(B)$ for all $A\subseteq B$, and that satisfies $ \nu(\varnothing)=0$.

The notion of capacity does not have a universally agreed definition. Our definition is consistent with the one used in \cite{MM04}, standard in decision theory. In some papers, a capacity is required to satisfy continuity from below and continuity from above on closed sets; see e.g., \cite{WK90}. The upper copula capacity, the main capacity of interest in this paper, satisfies both of these continuity requirements, as we will see below. 

A measure $C$ on $[0,1]^n$ is a \textit{copula} if all of its marginals are the Lebesgue measure on $[0,1]$,  meaning $C(p_i^{-1}(A)) = \ell(A)$ for all $i\in[n]$, $A\in\mathcal{B}$. 
Let $\mathbb{C}_n$ denote the set of all copulas on $[0,1]^n$.

\section{The upper copula capacity and elementary properties}
\label{sect:uppercopula}

The upper copula capacity $\uC_n$ and the lower copula capacity $\lC_n$ are defined by 
\begin{equation*}
\uC_n (A)=\sup_{C\in \mathbb C_n} C (A) \quad \text{and} \quad \lC_n(A) = \inf_{C\in \mathbb C_n} C (A),  \quad \quad \text{for every }A\in\mathcal{B}_n.
\end{equation*}
It is clear that 
  $
  \uC_n(A) + \lC_n(A^c) = 1 
  $ for any $A\in \mathcal{B}_n$, 
  and therefore 
  we will focus on the upper copula capacity $\uC_n$ in this section. Symmetric statements can be made for $\lC_n$, and they will be omitted here. We next present some properties of the upper copula capacity.

  \begin{proposition}
  \label{prop:capacity}
 Let $n\in\mathbb{N}$. The set function $\uC_n$ is a capacity on $([0,1]^n, \mathcal{B}_n)$. It moreover satisfies the following properties:
 \begin{enumerate}[label={\rm (\roman*)}, ref={\rm (\roman*)}]
     \item \label{item:property-1} $\uC_n([0,1]^n) = 1$ (it is normalized);
     \item \label{item:property-2} $\uC_n(A_1\cup A_2) \leq \uC_n(A_1) + \uC_n(A_2)$ for all $A_1,A_2\in\mathcal{B}_n$ (it is subadditive);
     \item \label{item:property-3} for any increasing sequence of measurable sets $\{A_m, m\in\mathbb{N}\}$, $\lim_{m\to\infty} \uC_n(A_m) = \uC_n(\bigcup_{m\in\mathbb{N}} A_m)$ (it is continuous from below);
     \item \label{item:property-4} for any decreasing sequence of \textbf{closed} sets $\{F_m, m\in\mathbb{N}\}$, $\lim_{m\to\infty} \uC_n(F_m) = \uC_n(\bigcap_{m\in\mathbb{N}} F_m)$ (it is continuous from above on closed sets). 
 \end{enumerate}
  \end{proposition}

A capacity $\nu$ is \textit{submodular} (or \textit{2-alternating}) when $\nu(A_1\cup A_2) + \nu(A_1\cap A_2) \leq \nu(A_1) + \nu(A_2)$ holds for all $A_1,A_2\in\mathcal{B}_n$. 
Many considerations on robustness require to know whether an upper probability is submodular; see for instance, \cite{B82}. 
The upper copula capacity  $\uC_n$ is not submodular. We illustrate this for $n=2$ in the following example; the case of general $n$ is similar.
\begin{example}\label{ex:not-submodular}
To see this for $n=2$, 
let $A=[0,1/2]\times [0,1]$ and $B=[0,1]\times [0,1/2]$. 
  It is clear that $\uC_2(A\cup B)=1$ since, notably, the countermonotonic copula yields $C(A\cup B)=1$.
  Moreover, $\uC_2(A)=\uC_2(B)=1/2$, and $\uC_2(A\cap B)=1/2$. Therefore, $\uC_2(A\cap B) +\uC_2(A\cup B) =3/2 > \uC_2(A)+\uC_2(B)$. 
\end{example}  
 
Another consideration is whether for every measurable set, some copula yields the same value as the upper copula capacity. This is answered in the following proposition and subsequent discussion.  
  
  \begin{proposition}
\label{prop:max-is-attained}
  If $A\in\mathcal{B}_n$ is closed, then the supremum defining $\uC_n$ is attained, and thus $\uC_n(A)=\max_{C\in \mathbb C_n} C(A)$.
  \end{proposition}

If $A$ is not closed, however, the supremum is not necessarily attained by a copula. For example, consider the set $A=\{\mathbf{u}\in[0,1]^2 : u_1>u_2\}$. No copula $C\in\mathbb{C}_2$ yields $C(A)=1$. Although, consider the sequence $\{C_m:m\in\mathbb{N}\}\subseteq \mathbb{C}_2$ where $C_m$ is the uniform distribution on the two-segment line $\{y=x-1/m: x\in[1/m,1)\}\cup \{y=1-x: x\in(0,1/m)\}$, for every $m\in\mathbb{N}$. 
For each $m\in\mathbb{N}$, $C_m(A) = 1-\frac{1}{m}$, which corresponds to the probability mass for the first segment of the line. We have $\sup_{C\in\mathbb{C}_2} C(A)\geq \sup_{m\in\mathbb{N}} C_m(A) = 1$, thus $\uC_2(A)=1$. Similar examples can be designed for the case of general $n$. 

Closedness also cannot be omitted for continuity from above, as specified in Proposition~\ref{prop:capacity}\ref{item:property-4}, and we here give an example showing that   $\uC_n$ is not continuous from above in a general sense. For $\alpha \in [0,1]$, define the set 
  $
M_\alpha = \{\mathbf{u}\in(0,1]^2:  (1-\alpha) u_1<u_2<u_1\}.
  $ 
  Clearly, $M_\alpha$ decreases to $M_0=\varnothing$ as $\alpha \to 0$, and hence,  $\uC_2(M_0) = 0$. 
  Meanwhile,  $\uC_2(M_\alpha)=1$ for every $\alpha \in (0,1)$, and therefore $\lim_{\alpha \to 0} \uC_2(M_{\alpha}) = 1$. An example for the case of general $n$ is obtained similarly by constructing $M_\alpha \times [0,1]^{n-2}$ in place of $M_\alpha$.

The Lebesgue measure $\ell_n$ on $[0,1]^n$ is  a copula itself, viz.~the independence copula; and hence, $\lC_n \leq \ell_n \leq \uC_n$.

\begin{proposition}
\label{prop:lines}
    Let $A\in\mathcal{B}_n$ be closed and $n\geq 2$, then $\argmax_{C\in\mathbb{C}_n}C(A)$ contains a copula $C^*$ such that $C^* \perp \ell_n$. 
\end{proposition}

Proposition~\ref{prop:lines} shows that we could replace $\mathbb{C}_n$ by the set of copulas singular with respect to the Lebesgue measure in corresponding dimension; that is, we do not need to check every copula, only singular ones. Proposition~\ref{prop:lines} assumes that $A$ is closed, but, for a general Borel $A$, the mass of $A$ under each copula can be approximated from below by the mass on closed subsets of $A$. 

There is a dual representation of $\uC_n$, which is not useful in computing $\uC_n(A)$ but helpful in verifying some of its properties. 
\citet[Corollary 2.18]{K84} showed that for every $A\in \mathcal B_n$,
\begin{equation}
    \label{eq:duality} 
\uC_n(A)= \inf
\left\{\sum_{i=1}^n \int_0^1 f_i(x)\d x :\bigoplus_{i=1}^n f_i \ge \id_A;~f_1,\dots,f_n\in \mathcal{F}_{\rm BM}\right\}.
\end{equation}
where $\bigoplus_{i=1}^n f_i$ is the mapping 
$(x_1,\dots,x_n)\mapsto 
 \sum_{i=1}^n f_i(x_i) $  on $[0,1]^n$
 and $\mathcal{F}_{\rm BM}$ is the set of all real-valued bounded measurable functions on $[0,1]$. 

\section{Geometry of events under the upper copula capacity}
\label{sect:geometry}

We next study the measure of events with respect to the upper copula capacity. This is done by examining how some geometrical characteristics of events $A\in\mathcal{B}_n$ translate to $\uC_n(A)$. Leveraging these properties then allows to compute $\uC_n(A)$. We first give elementary results in Section~\ref{subsect:geo-elementary}. In Section~\ref{subsect:geo-kellerer}, we examine the two-dimensional problem and discuss its connection to the marginal problem of \cite{S65}. In Section~\ref{subsect:geo-convex}, we give results for convex sets; in Section~\ref{subsect:geo-decreasing}, we give results for decreasing sets. 
In all of these sections, we only discuss $\uC_n$; results for the lower copula capacity follow by $\lC_n(A) = 1-\uC_n(A^{\rm c})$ for any $A\in\mathcal{B}_n$. 

\subsection{Elementary geometrical considerations}
\label{subsect:geo-elementary}

The upper copula capacity of some sets follow directly from the uniform marginals.  

\begin{proposition}
\label{prop:intervals} 
      The relation $\uC_{k+m}(A\times B) = \min(\uC_k(A), \uC_m(B))$ holds for any $A\in\mathcal{B}_k$ and $B\in\mathcal{B}_m$, $k,m\in\mathcal{\mathbb{N}}$. 
Consequently, 
\begin{equation*}
    \uC_n\left(\prod_{i=1}^n A_i\right) = \min_{i\in[n]} \ell(A_i)\quad \text{for any }  A_1,\ldots,A_n\in\mathcal{B}.
\end{equation*} 
\end{proposition}
Some of our results below hold only for dimension $n=2$, but they may nonetheless prove useful in higher dimensions by referring to Proposition~\ref{prop:intervals}. For example, if $A=A_1\times A_2$ for some $A_1,A_2\in\mathcal{B}_2$, then 
$\uC_4(A) = \min(\uC_2(A_1),\uC_2(A_2))$ and the two-dimensional results will be applicable.

We next examine a broad class of transformations for a set that leaves its measure under the upper copula capacity unchanged. 
Let $\Phi_U$ be the set of rearrangements of the identity, i.e.~the set of measurable functions $\varphi$ such that $\varphi(U)$ follows a uniform distribution on $(0,1)$ whenever $U$ does. Such functions are notably studied in \cite{PW15}. 
In our next result, we show that $\uC_n$ is invariant to any rearrangements of the identity applied componentwise.

For measurable functions $\varphi_1,\ldots,\varphi_n\in\Phi_U$,
let the function $\varphi:[0,1] ^n \to [0,1]^n $ be given by
$(u_1,\dots,u_n)\mapsto (\varphi_1(u_1),\dots,\varphi_n(u_n)).$

\begin{theorem}
\label{th:rearrangements}
For any set $A\in\mathcal{B}_n$ and  $\varphi_1,\ldots,\varphi_n\in\Phi_U$, we have $\uC_n(A) = \uC_n(\varphi^{-1}(A))$.
\end{theorem}

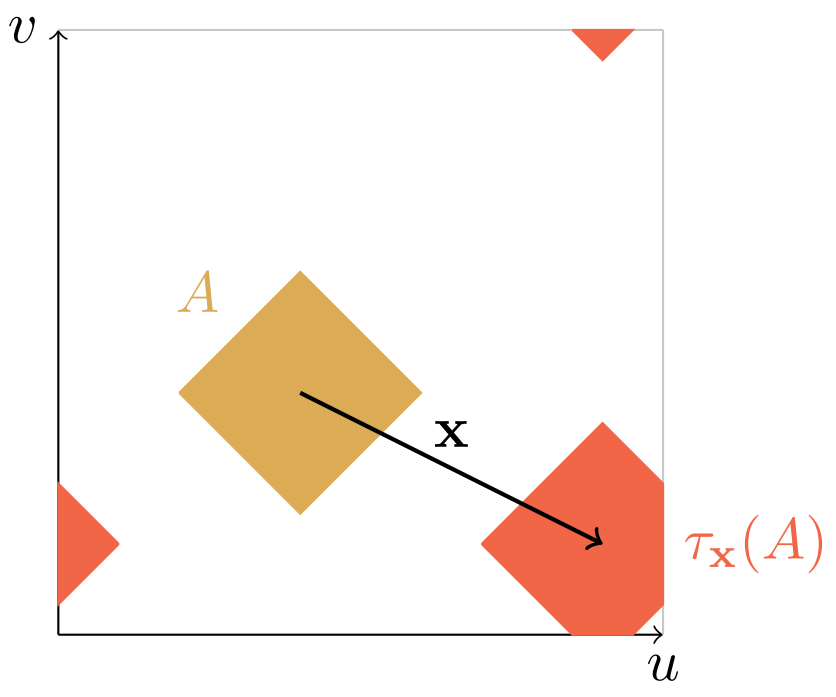
\begin{figure}[tb]
\centering
      \begin{tikzpicture}[scale = 4]
          \draw[->] (0,0) -- (1,0) node[below]{$u$} ;
		\draw[->] (0,0) -- (0,1) node[left]{$v$} ;
        \draw[Black!25] (1,0) -- (1,1);
        \draw[Black!25] (0,1) -- (1,1);

        \filldraw[Dandelion!85!black] (0.2,0.4) -- node[midway, above left]{$A$} (0.4,0.6) -- (0.6,0.4) -- (0.4,0.2) -- (0.2,0.4); 
        \filldraw[Red!75] (0.7,0.15) -- (0.9,0.35) -- (1,0.25) -- node[midway, right]{$\tau_{\mathbf{x}}(A)$} (1,0.05) -- (0.95,0) -- (0.85,0) -- (0.7,0.15); 
        \filldraw[Red!75] (0.1,0.15) -- (0,0.05) -- (0,0.25) -- (0.1,0.15);
        \filldraw[Red!75] (0.85,1) -- (0.9,0.95) -- (0.95,1) -- (0.85,1);
        \draw[->, thick] (0.4,0.4) -- node[midway, above]{$\mathbf{x}$} (0.9,0.15);
    \end{tikzpicture}
\caption{An example of translation performed on some set $A$ in yellow. The set $\tau_{\mathbf{x}}(A)$ is depicted in red; we have $\uC_2(\tau_{\mathbf{x}}(A)) = \uC_2(A)$.}
\label{fig:trans}
\end{figure}

An example of transformations on $A$ realizable as the preimage of componentwise rearrangement are translations where we
suppose the space loops back on itself in every coordinate, as if $[0,1]^n$ were the surface of a $n$-dimensional torus. Figure~\ref{fig:trans} provides an example.  
More formally, we define the translation operator $\tau_{\mathbf{x}}: \mathcal{B}_n \mapsto \mathcal{B}_n$ by $\tau_{\mathbf{x}}(A) = \left\{\langle \mathbf u + \mathbf x \rangle: \mathbf{u}\in A\right \}$
 for any $\mathbf{x} \in \mathbb{R}^n$, where $\langle \mathbf z\rangle$ denotes the fractional part of $\mathbf z\in\mathbb{R}^ n$, applied componentwise.  

 \begin{corollary}
\label{prop:translations}
     We have 
$\uC_n(\tau_{\mathbf{x}}(A)) = \uC_n(A)$   for any $A\in\mathcal{B}_n$ and $\mathbf{x}\in\mathbb{R}^n$.
 \end{corollary}

Another transformation to be considered is taking the closure of a set.  
For any $A\in\mathcal{B}_n$, let  $A^{\circ}$ denote its interior
and $\overline{A}$ its closure. 
We say that a set $A$ is \textit{strictly star-shaped}  if there exists $\mathbf a\in A^{\circ}$ such that $t \overline{A} + (1-t)\mathbf a\subseteq A^{\circ}$ for every $t\in(0,1)$. Equivalently,
for any $\mathbf a'\in \overline A$, the open line segment from $\mathbf a$ to $\mathbf a'$ is contained in $A^{\circ}$.
\begin{proposition}
\label{lem:closure} 
If  $A\in\mathcal{B}_n$ is strictly star-shaped, then
\begin{equation*}
\uC_n(A) = \uC_n(\overline{A}) = \uC_n(A^{\circ}) = \uC_n(\overline{A^{\circ}}).
\end{equation*}
\end{proposition}

    All convex sets with nonempty interior are notably strictly star-shaped. 
    This condition
    cannot be deleted in Proposition~\ref{lem:closure}. For instance, take $A =([0,1]\cap\mathbb{Q})\times [0,1]$, where $\mathbb{Q}$ denotes the rational numbers. By Proposition~\ref{prop:intervals}, we have $\uC_2(A) = 0$, whereas $\uC_2(\overline{A}) = 1$ since $\overline{A}=[0,1]^2$. Note that in this case  $A^{\circ} = \varnothing$; thus, $A$ is not strictly star-shaped.

\subsection{Two-dimensional sets and the marginal problem}
\label{subsect:geo-kellerer}

For a closed $A\in\mathcal{B}_2$, verifying whether $\uC_2(A) = 1$ amounts to determining whether there exists a copula whose support is contained in $A$. This problem, for general marginals, not just uniform, has been answered in two dimensions by \cite{S65}, \cite{E78}, \cite{GR84} and \cite{S93}. Results by \cite{K84} also solve this problem; in fact, the following proposition is directly derived from his results.

\begin{proposition}
\label{prop:dim2}
\label{prop:strip-covers}
For any closed $A\in\mathcal{B}_2$,
\begin{align}
\uC_2(A)
 &= \inf\{\ell(S)+\ell(T): A\subseteq(p_1^{-1}(S)\cup p_2^{-1}(T)),\; S,T\in\mathcal{B}\}
 \label{eq:strip-covers}\\
 &= 1-\sup_{\substack{K\subseteq[0,1]\\K \text{\rm \,is closed}}}\{\ell(K)-\ell(\pi_A(K))\},
 \label{eq:hall-deficiency}
\end{align}
where $\pi_A(K)=\{y:(x,y)\in A\text{ for some }x\in K\}$. Therefore, a copula concentrated on $A$ exists if and only if $\ell(K)\leq\ell(\pi_A(K))$ for every closed $K\subseteq[0,1]$.
\end{proposition}

Note that 
$\ell(S)+\ell(T)=\int \id_{S } \d \ell  
+ \int \id_{ T} \d \ell  $
and  $A\subseteq(p_1^{-1}(S)\cup p_2^{-1}(T))$ 
is equivalent to
 $\id_{S } (s) +\id_{  T}(t) \ge \id_A(s,t)$ for all $s,t\in[0,1]$. 
Comparing \eqref{eq:duality}
with \eqref{eq:strip-covers},
the functions $f_1$ and $f_2$ in \eqref{eq:duality} are replaced by the indicator functions, which makes the infimum weakly bigger; yet equality holds. This is  true only in dimension $2$.   For $n\geq 3$, \eqref{eq:duality} also gives the following upper bound:
 For any closed $A\in\mathcal{B}_n$,
 \begin{equation}
     \uC_n(A) \leq \inf\left\{\sum_{i=1}^n \ell(S_i) : A\subseteq \bigcup_{i=1}^n p_i^{-1}(S_i), \; S_1,\ldots,S_n\in\mathcal{B} \right\} .
     \label{eq:kellerer-dim3}
 \end{equation}
The inequaltiy in \eqref{eq:kellerer-dim3} is famously not an equality in general for the case $n\geq 3$; see \citet[Remark~5]{GR84} and \citet[Example~3.4]{K84}.

In \eqref{eq:hall-deficiency}, $\pi_A(K)$ can be replaced by $\{x : (x,y) \in A \;\text{for some }y\in K \}$ by symmetry. 
By Proposition~\ref{lem:closure}, the identities in \eqref{eq:strip-covers} and \eqref{eq:hall-deficiency} moreover holds for all $A$ that are strictly star-shaped. 

The criteria in Proposition~\ref{prop:dim2} are practical for applications as one can deploy geometrical arguments to argue that $\ell(K) \leq \ell(\pi_A(K))$ for all closed $K$. Other criteria from the results on the marginal problems would similarly extend to $\uC_n$; however, they fall beyond our geometrical analysis, as for example some require verification for all pairs of continuous bounded functions and this is essentially the same level of complexity as verifying each copula.


\begin{example}
\label{ex:projection-deficit}
Consider
$
A=([0,3/4]\times[0,1/4])\cup([3/4,1]\times[1/4,1]).
$
Both coordinate projections $p_1(A)$ and $p_2(A)$ are the whole interval $[0,1]$. Nevertheless, a cover consisting of $S=[3/4,1]$ and $T=[0,1/4]$ yields $\ell(S) + \ell(T) = 1/2$.  Hence, $\uC_2(A) \leq 1/2$ by \eqref{eq:strip-covers}. The comonotonic copula $C^{\rm co}$ gives $C^{\rm co}(A) = 1/2$; thus, $\uC_2(A)\geq 1/2$. We conclude that $\uC_2(A) = 1/2$. 
\end{example}

\subsection{Convex sets}
\label{subsect:geo-convex}

For any nonempty set $A\in\mathcal{B}_n$ and $i\in[n]$, define $a_i(A) = \inf p_i(A)$ and $b_i(A) = \sup p_i(A)$. Write $I_A =\prod_{i=1}^n [a_i(A),b_i(A)] $, and let $I_A=\varnothing$ if $A=\varnothing$.

  \begin{theorem}
  \label{th:convex}
    For any convex set $A\in\mathcal{B}_2$, 
\begin{equation}
    \uC_2(A) = \uC_2(I_A) = \min_{i\in\{1,2\}}(b_i(A)-a_i(A)).
    \label{eq:convex}
\end{equation}
 \end{theorem}

The proof of Theorem~\ref{th:convex} relies on the statement that any closed convex set touching the four edges of $[0,1]^2$ contains the support of a copula. We provide, in Appendix~\ref{sect:construction-copula}, an explicit construction of such a copula for any closed convex set touching the four edges. 


Theorem~\ref{th:convex}
implies that $\uC_2(A)=1$
for any convex $A$ that touches all edges of $[0,1]^2$. This is related to the plank-covering problem studied by \cite{B51}.  By combining \eqref{eq:strip-covers} and \eqref{eq:convex}, we obtain the following geometric result: \textit{The sum of the relative widths of strips required to cover a two-dimensional convex set of nonzero area is at least one, where strips can be in two directions only and each of their width is measured as a fraction of the set’s width in this direction.} Taking relative widths means that the convex set is rescaled to have width 1 in the two directions. 
This geometric statement is the (proved)  two-directional case of \cite{B51}'s plank conjecture in two dimensions with   planks along the axes. Our proof of Theorem~\ref{th:convex} is thus, combined with Proposition~\ref{prop:dim2},  an alternative proof to this case of the conjecture. See \citet[Theorem 7.1]{AKP19} for another proof of this specific case relying on a completely different approach.

\begin{figure}[tb]
	\centering
	\begin{tikzpicture}[scale = 3]

    \draw[Cyan!30, very thick, ->] (0.6,0.9,0.5) to [out=0, in=180] (2,0.5);
             \draw[Cyan!30, very thick, ->] (0.6,0.1,0.5) to [out=0, in=180] (2,0.5);
              \draw[Blue!30, very thick, ->] (0.7,0.8,0.5) to [out=0, in=180] (2,0.5);
               \draw[Blue!30, very thick, ->] (0.7,0.2,0.5) to [out=0, in=180] (2,0.5);
	  \draw[Cyan] (0.4,0.1,0.5) -- (0.5,0.1,0.6) -- (0.6,0.1,0.5) -- (0.5,0.1,0.4) -- (0.4,0.1,0.5);
            \draw[Cyan] (0.4,0.9,0.5) -- (0.5,0.9,0.6) -- (0.6,0.9,0.5) -- (0.5,0.9,0.4) -- (0.4,0.9,0.5);
             \draw[Blue] (0.3,0.2,0.5) -- (0.5,0.2,0.7) -- (0.7,0.2,0.5) -- (0.5,0.2,0.3) -- (0.3,0.2,0.5);
            \draw[Blue] (0.3,0.8,0.5) -- (0.5,0.8,0.7) -- (0.7,0.8,0.5) -- (0.5,0.8,0.3) -- (0.3,0.8,0.5);	
        \draw[->] (0,0,0) -- (0,0,1) node[above left]{$u_{1}$} ;
		\draw[->] (0,0,0) -- (0,1,0) node[left]{$u_{3}$} ;
		\draw[->] (0,0,0) -- (1,0,0) node[below]{$u_{2}$} ;
		\filldraw[Sepia, opacity=0.17] (0,0.5,0.5) -- (0.5,0.5,1) -- (0.5,1,0.5);
		\filldraw[Sepia, opacity=0.05] (0.5,1,0.5) -- (1,0.5,0.5) -- (0.5,0.5,1);
		\filldraw[Sepia, opacity=0.10] (0.5,0,0.5) -- (1,0.5,0.5) -- (0.5,0.5,1);
		\filldraw[Sepia, opacity=0.25] (0.5,0,0.5) -- (0,0.5,0.5) -- (0.5,0.5,1);
            \filldraw[Sepia, opacity=0.35] (0,0.5,0.5) -- (0.5,0.5,0) -- (0.5,1,0.5);
		\filldraw[Sepia, opacity=0.35] (0.5,0.5,0) -- (0.5,1,0.5) -- (1,0.5,0.5);
		\filldraw[Sepia, opacity=0.30] (0.5,0.5,0) -- (1,0.5,0.5) -- (0.5,0,0.5);
		\filldraw[Sepia, opacity=0.40] (0.5,0.5,0) -- (0,0.5,0.5) -- (0.5,0,0.5);
        
            \draw[Sepia] (0,0.5,0.5) -- (0.5,0.5,1);
            \draw[Sepia] (0.5,0.5,1) -- (0.5,1,0.5);
		\draw[Sepia] (0.5,1,0.5) -- (1,0.5,0.5) ;
           \draw[Sepia] (1,0.5,0.5) -- (0.5,0.5,1);
		\draw[Sepia] (0.5,0,0.5) -- (1,0.5,0.5);
       		\draw[Sepia] (0.5,0,0.5) -- (0.5,0.5,1);
		\draw[Sepia] (0.5,0,0.5) -- (0,0.5,0.5);
       \draw[Sepia] (0,0.5,0.5) -- (0.5,0.5,1);
            \draw[Sepia] (0.5,0.5,1) -- (0,0.5,0.5) -- (0.5,1,0.5);
            \draw[Sepia] (0.5,0.5,1) -- (0.5,0,0.5);

          \draw[->, xshift = 2.1cm] (0,0) -- (1,0) node[below]{$u_{1}$} ;
		\draw[->, xshift = 2.1cm] (0,0) -- (0,1) node[left]{$u_{2}$} ;
        \draw[Black!25, xshift=2.1cm] (1,0) -- (1,1);
        \draw[Black!25, xshift=2.1cm] (0,1) -- (1,1);
        \filldraw[draw=Blue, fill = Blue!30, xshift = 2.1cm] (0.3,0.5) -- (0.5,0.7) -- (0.7,0.5) -- (0.5,0.3) -- (0.3,0.5); 
                  \filldraw[draw=Cyan, fill = Cyan!30, xshift = 2.1cm](0.4,0.5) -- (0.5,0.6) -- (0.6,0.5) -- (0.5,0.4) -- (0.4,0.5); 
        \draw[dashed, Cyan!30, xshift=2.1cm] (0,0.4) -- (1,0.4);
        \draw[dashed, Cyan!30, xshift=2.1cm] (0,0.6) -- (1,0.6);
        \draw[dashed, Cyan!30, xshift=2.1cm] (0.4,0) -- (0.4,1);
        \draw[dashed, Cyan!30, xshift=2.1cm] (0.6,0) -- (0.6,1);
        \draw[dashed, Blue!30, xshift=2.1cm] (0,0.3) -- (1,0.3);
        \draw[dashed, Blue!30, xshift=2.1cm] (0,0.7) -- (1,0.7);
        \draw[dashed, Blue!30, xshift=2.1cm] (0.3,0) -- (0.3,1);
        \draw[dashed, Blue!30, xshift=2.1cm] (0.7,0) -- (0.7,1);

        \draw[Cyan, xshift = 2.1cm] (0.4, -0.08) -- (0.4, -0.1) -- node[below] {0.2} (0.6,-0.1) -- (0.6, -0.08);
        \draw[Blue, xshift = 2.1cm] (0.3, -0.28) -- (0.3, -0.3) -- node[below] {0.4} (0.7,-0.3) -- (0.7, -0.28);
        
        \draw[Cyan, xshift = 2.1cm] (1.08, 0.4) -- (1.1, 0.4) -- node[right] {0.2} (1.1,0.6) -- (1.08, 0.6);
        \draw[Blue, xshift = 2.2cm] (1.28, 0.3) -- (1.3, 0.3) -- node[right] {0.4} (1.3,0.7) -- (1.28, 0.7);

	\end{tikzpicture}
	\caption{(Left) The regular octahedron $A$ from Example~\ref{ex:octahedron}. (Right) Sets $\{|U_1-0.5|+|U_2-0.5|\leq 0.1\}$, in light blue, and $\{|U_1-0.5|+|U_2-0.5|\leq 0.2\}$, in violet. It is impossible to attribute 0.2 units of probability mass to the light blue region and 0.4 to the violet region while respecting both uniform margins.}
	\label{fig:octahedron}
\end{figure}
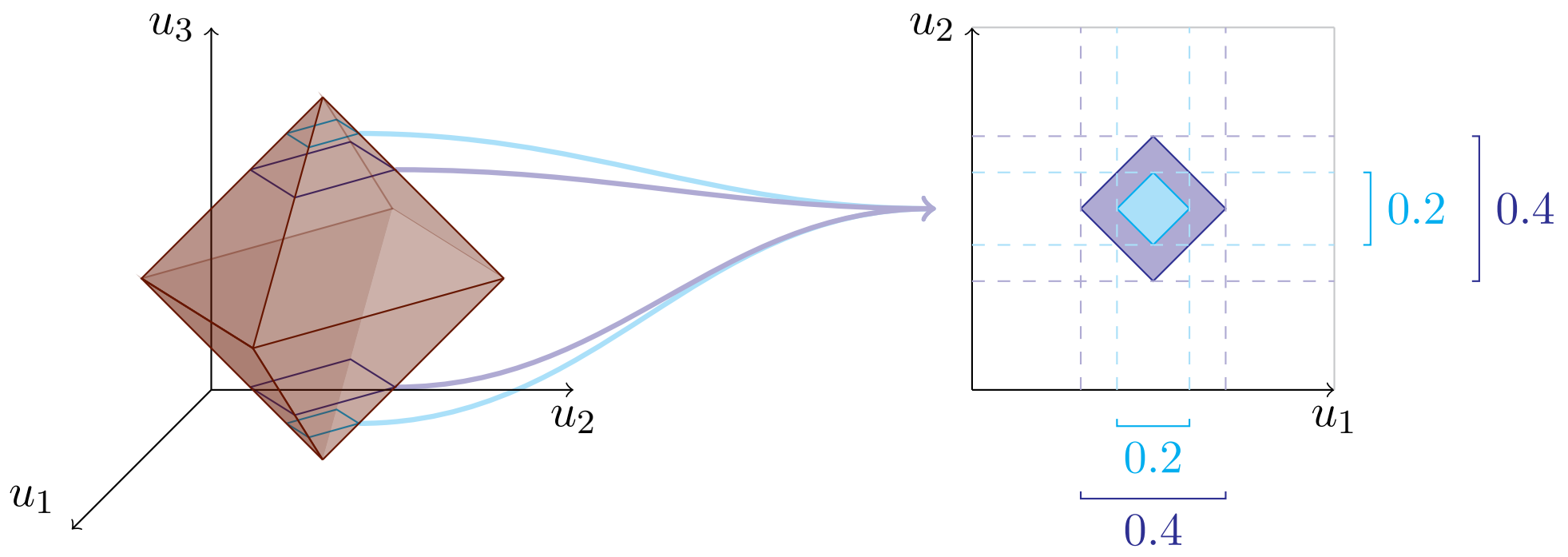

The  conclusion that $\uC_n(A)=1$ for convex $A$ touching all sides of $[0,1]^n$  does not hold for $n\ge 3$. Indeed, $\uC_n(\Delta_n)=2/n$ for all $n\ge 2$,
where $\Delta_n$ is the dimensional simplex in $\R^n$; see \citet[Section 3]{AHKM02} formulated in the context of fractional plank covering.  
We provide another example below to illustrate the interesting structure of $\uC_n$ in dimension $n=3$.
\begin{example}
\label{ex:octahedron}
Let $A$ be the regular octahedron inscribed in $[0,1]^3$ such that each of its vertices touches one face of the cube. This set is convex. If \eqref{eq:convex} were to extend to dimension $n=3$, we would have $\uC_3(A) = 1$ since $I_A=[0,1]^3$. We will argue that  $\uC_3(A) < 1$. Suppose that there exists $C^*\in\mathbb{C}_3$ such that $C^*(A) = 1$, and let $(U_1,U_2,U_3)\sim C^*$. From the uniform distribution of $U_3$, one requires $C^*(\{U_3 \geq 0.9\}\cup \{U_3 \leq 0.1\}) = 0.2.$
On $\{U_3 \geq 0.9\}\cup \{U_3 \leq 0.1\}$, we have $|U_1-0.5|+|U_2-0.5|\leq 0.1$ almost surely given the octahedron constraint; thus this event has probability at least 0.2. It also has probability at most 0.2, since it is contained in $\{|U_1-0.5|\le 0.1\}$.
Hence, 
\begin{equation}
C^*(|U_1-0.5|+|U_2-0.5|\leq 0.1) = 0.2.    
\label{eq:octahedron-1}
\end{equation}
 Similarly, one requires $C^*(\{U_3 \geq 0.8\}\cup \{U_3 \leq 0.2\}) = 0.4$, 
thus 
\begin{equation}
    C^*(|U_1-0.5|+|U_2-0.5|\leq 0.2) = 0.4. 
    \label{eq:octahedron-2}
\end{equation}
Individually, each of \eqref{eq:octahedron-1} and \eqref{eq:octahedron-2} are valid requirements on $C^*$; however, together, they violate the uniformness of $U_1$ and $U_2$ as any probability mass on $\{0.1< |U_1-0.5|+|U_2-0.5|\leq 0.2\}$ either adds to $\{|U_1-0.5|\leq 0.1\}$ or $\{|U_2-0.5|\leq 0.1\}$, the two sets already bearing $0.2$ probability mass from the requirement \eqref{eq:octahedron-1}, the maximum permitted by the uniform margins. See Figure~\ref{fig:octahedron} for an illustrative support to this argument.
We are left with $0.4-0.2 = 0.2$ unattributable probability mass, meaning $C^*$ is unachievable. Consequently, $\uC_3(A)< 1$. We derive the exact value of $\uC_3(A)$ in Example~\ref{ex:inscribed-octahedron-exact} below. 
\end{example}

In the next subsection, we show how to combine Theorem~\ref{th:rearrangements} and exact-computation results for a class of decreasing sets to obtain exact computation of some multidimensional convex sets.

  \subsection{Decreasing sets}
\label{subsect:geo-decreasing}
  
  We say that a set $A\subseteq [0,1]^n$ is a decreasing set if
$\mathbf x \in A $ implies $\mathbf y\in A $ for all
$ \mathbf y\in [0, 1]^n$ with $\mathbf y
\le  \mathbf x$ (componentwise).
For some applications, the sets needing to be measured by the upper-copula capacity are limited to decreasing sets, see Section~\ref{sect:app}. 

The next result establishes homogeneity of  $\uC_n$ for decreasing sets.
  \begin{proposition}
  \label{th:decreasing-sets}
  For any decreasing Borel set $A\subseteq [0,1]^n$ and $\lambda \in [0,1]$, 
  we have 
  $\uC_n(\lambda A)=\lambda \uC_n(  A)$.
  \end{proposition}

A specific type of decreasing set is the simplex obtained as the convex hull of points on each axis and the origin. This set may be interpreted as feasible allocations of resources.  
By the joint-mixability theorem of \citet[Theorem~3.1]{WW16}, uniform distributions on intervals of respective lengths $\lambda_1, \ldots,\lambda_n$ admit a coupling with constant sum if and only if $\max_{i\in[n]} \lambda_i\leq\frac{1}{2}\sum_{i=1}^n \lambda_i$. We employ this result to derive the capacity of the simplex (clipped to $[0,1]^n$). 


\begin{theorem}
\label{th:weighted-uniform-sums}
For any $\boldsymbol{\lambda} = (\lambda_1,\ldots, \lambda_n)\in(0,\infty)^n$ and $s\geq0$,
\begin{equation*}
\uC_n\left(\left\{\mathbf{u}\in[0,1]^n:\sum_{i=1}^{n} \lambda_i u_i\leq s\right\}\right)
=\min\left\{\frac{s}{L(\boldsymbol{\lambda})},1\right\}, 
\end{equation*}
where
$
L(\boldsymbol{\lambda})=\max(\max_{i\in[n]} \lambda_i,\frac{1}{2}\sum_{i=1}^n \lambda_i). 
$
\end{theorem}

The noteworthiness of Theorem~\ref{th:weighted-uniform-sums} lies in that it provides exact computation, for simplices, in dimension $n\geq 3$; for the case $n=2$, one could instead refer to Theorem~\ref{th:convex}. Combined with the transformations from Section~\ref{subsect:geo-elementary}, Theorem~\ref{th:weighted-uniform-sums} thus permits to compute the upper copula capacity of many sets in dimension $n\geq 3$, as we illustrate in the following example. 


\begin{example}
\label{ex:inscribed-octahedron-exact}
Consider the regular octahedron from Example~\ref{ex:octahedron} (the one in Figure~\ref{fig:octahedron}). It can be expressed as $A=\{\mathbf{u}\in[0,1]^3: \sum_{i=1}^3 |u_i - 1/2|\leq 1/2 \}$. By considering the rearrangement $\varphi\in\Phi_U$ given by $\varphi(v) = 2|v-1/2| $, we have
\begin{equation*}
    \uC_3(A) = \uC_3\left(\left\{u\in[0,1]^3: \sum_{i=1}^3 \varphi(u_i) \leq 1 \right\}\right) = \uC_3\left(\left\{u\in[0,1]^3: \sum_{i=1}^3 u_i \leq 1 \right\}\right),
\end{equation*}
 by Theorem~\ref{th:rearrangements}. 
Thus, we obtain $\uC_3(A) = 2/3$ by applying Theorem~\ref{th:weighted-uniform-sums} with $\lambda_1 = \lambda_2 = \lambda_3 = s = 1$. 
\end{example}

\section{Integration with respect to the upper/lower copula capacity}
\label{sect:INT}


The Choquet integral of a measurable function $\psi:[0,1]^n\to \R$ with respect to a normalized capacity $\nu$ on $([0,1]^n,\B_n)$ is defined as 
\begin{equation}
\int \psi \d \nu = \int_0^\infty \nu(\psi>z) \d z
+ \int_{-\infty}^0\left( \nu(\psi>z)-1\right) \d z 
 = \int_{\mathbb{R}} \left( \nu(\psi>z) - \id_{\{z<0\}} \right) \d z,
\label{eq:def-choquet}
\end{equation}
assuming that the above term is well defined; see, e.g., \citet[Definition 4.1]{MM04}. Let $\uJ_n$ and $\lJ_n$ denote integration with respect to the upper and lower copula capacity respectively, meaning $\uJ_n(\psi) = \int \psi \mathrm{d}\uC_n$ and $\lJ_n(\psi) = \int \psi \mathrm{d}\lC_n$ for any measurable function $\psi$ in the appropriate dimension. The two functionals are connected through the relation   $\lJ_n(\psi) = -\uJ_n(-\psi)$; and hence, for any measurable function $\psi$ that is bounded on $[0,1]^n$,
\begin{equation}
       \lJ_n(\psi) = m - \uJ_n(m-\psi)\quad \text{for every }m\in\mathbb{R}.
\label{eq:link-lJ-uJ}
\end{equation}

\subsection{Integration with the upper copula capacity vs.~worst-case integration}
\label{subsect:worst-case}

Let us also define $\uK_n$ and $\lK_n$ respectively as the worst-case and best-case integrals of $\psi$, meaning $\uK_n(\psi) = \sup_{C\in\mathbb{C}_n} \int \psi \mathrm{d}C$ and $\lK_n(\psi) = \inf_{C\in\mathbb{C}_n} \int \psi \mathrm{d} C$ for any measurable function $\psi$ in the appropriate dimension. It is straightforward to see from \eqref{eq:def-choquet} that 
\begin{equation}
    \lJ_n \leq \lK_n \leq \uK_n \leq \uJ_n,
    \label{eq:layers}
\end{equation}
whenever all four quantities are well-defined. 
 The inequality $\uK_n(\psi) \leq \uJ_n(\psi)$ is not an equality in general; idem for $\lJ_n \leq \lK_n$. To justify this, suppose that $\uK_n = \uJ_n$. Since 
 $\uK_n$
 is the supremum of linear functions,
 this mapping is subadditive. 
 For a capacity $\nu$, 
the  mapping  $\psi \mapsto \int \psi \d \nu$ is subadditive if and only if $\nu$ is submodular; see \citet[Section 54]{C54}. 
This yields a contradiction, as we have seen in Example \ref{ex:not-submodular} that $\uC_n$ is not submodular.   
The difference between $\uJ_n$ and $\uK_n$ is examined in the following proposition. 

\begin{proposition}
\label{prop:JK}
    For a bounded, measurable and upper semicontinuous $\psi:[0,1]^n\to \mathbb{R}$, the following are equivalent
    \begin{enumerate}[label={\rm (\roman*)}, ref={\rm (\roman*)}]
        \item \label{item:JK-eq} $\uJ_n(\psi) = \uK_n(\psi)$;
        \item \label{item:JK-att} there exists a copula $C^*\in\mathbb{C}_n$ such that $C^*(\psi> z) = \uC_n(\psi>z)$ for every $z\in\mathbb{R}$. 
    \end{enumerate}
\end{proposition}

Functions $\psi$ satisfying the criterion in \ref{item:JK-att} are easy to define by first selecting a copula $C^*$ and then making sure that every superlevel set $\{\mathbf{u}\in[0,1]^n : \psi(\mathbf{u}) > z\}$ yields the same measure under $C^*$ and $\uC_n$. Conversely, if, for a given $\psi$, one can show that 
any copula $C^{\dagger}$ attaining the supremum for  $\uC_n(\psi>z_1)$, at some $z_1\in\mathbb{R}$, has $C^{\dagger}(\psi>z_2) < \uC_n(\psi>z_2)$ for some $z_2\in\mathbb{R}$, then the strict inequality $\uJ_n(\psi)> \uK_n(\psi)$ holds by Propositions~\ref{prop:JK} and \eqref{eq:layers}.
Indicator functions for closed set (thus the sup is attained by   Proposition~\ref{prop:max-is-attained}) always satisfy the criterion in \ref{item:JK-att}, so equality always hold for them. 

In decision theory, 
the functions 
$\lJ_n$, $\uJ_n$, $\lK_n$,
 and $\uK_n$ correspond to classical models for decision making under ambiguity
consistent with a   lack of knowledge regarding dependence. 
  The functional $\lJ_n$  corresponds  to Choquet expected utility \citep{S89} under capacity $\lC_n$. To be precise, for a utility function $\mathcal U: \X\to\R$, where $\X$ is a fixed set of consequences, and two random payoffs $X,Y:[0,1]^n \to \X$,
  the Choquet expected utility model is given by
  $$
X\succsim Y ~\iff~\int \mathcal U (X) \d \lC_n \ge   \int \mathcal U (Y) \d \lC_n ~\iff~    \lJ_n  (\mathcal U  (X))   \ge   \lJ_n (\mathcal U  (Y)).
  $$
  A concrete example is, with $\X=\R$,  $X(u_1,\dots,u_n)=F_1^{-1}(u_1)+\dots+F_n^{-1}(u_n)$
which corresponds to the total portfolio payoff  with  marginals $F_1,\dots,F_n$,  
and the set $\mathbb C_n$ describes the ambiguity on dependence.
On the other hand, $\lK_n$   corresponds  to maxmin expected utility \citep{GS89} in which the set of prior is   $\mathbb{C}_n$.  To be precise, under the same setting above, the maxmin expected utility model is given by
  $$
X\succsim Y ~\iff~\inf_{C\in \mathbb C_n} \int \mathcal U (X) \d C \ge \inf_{C\in \mathbb C_n}   \int \mathcal U (Y) \d C ~\iff~    \lK_n  (\mathcal U  (X))   \ge   \lK_n (\mathcal U  (Y)).
  $$
Both \cite{S89}'s and \cite{GS89}'s model are classic for decision making under ambiguity and have solid decision theoretic foundations. Using the sign convention of \cite{DDGKL06} and \cite{MFE15} in risk management, $\lJ_n $ and $\lK_n $ 
become $\uJ_n $ and $\uK_n $.
Proposition~\ref{prop:JK} thus indicates precisely for which $\mathcal U(X)$, satisfying the assumptions, the two models agree.

\subsection{Integrating quasiconcave functions}

A function $\psi:[0,1]^n\to \mathbb{R}$ is \textit{quasiconcave} if its superlevel sets $\{\mathbf{u}\in[0,1]^n : \psi(\mathbf{u}) > z \}$ are convex for every $z\in\mathbb{R}$. In dimension $n=2$, quasiconcave $\psi$ have an exact integration formula with respect to the upper copula capacity, provided in the following proposition. 

\begin{proposition}
\label{prop:integr-gravitational}
   Let $\psi:[0,1]^2\to\mathbb{R}$ be measurable and quasiconcave. We have
    \begin{equation*}
        \uJ_2(\psi) = \inf\psi + \int_{\inf \psi}^{\sup \psi}\min_{i\in\{1,2\}}(\ukappa_i(z) - \lkappa_i(z)) \, \mathrm{d}z,
    \end{equation*}
    where 
    \begin{equation*}
        \ukappa_i(z) = \sup\left\{u_i : \psi(\mathbf{u})>z, \mathbf{u}\in[0,1]^2\right\} \quad\text{and}\quad \lkappa_i(z) = \inf\left\{u_i : \psi(\mathbf{u})>z, \mathbf{u}\in[0,1]^2\right\}
    \end{equation*}
    for $i\in\{1,2\}$ and every $z<\sup\psi$.  
\end{proposition}


Similarly, one may integrate quasiconvex bounded funcions with respect to the lower copula capacity. 
In the following example, we illustrate an application of Proposition~\ref{prop:integr-gravitational}.

\begin{figure}[ht]
    \centering
    \begin{tikzpicture}

    \begin{axis}[ xmin=-0.1, xmax = 1, zmin=0, zmax=1, hide axis, scale =0.8]
\draw[->, very thick] (0,0,0.5) -- (0,1,0.5)node[left]{$\boldsymbol{u}_2$};
\addplot3 [
domain=0:1,
domain y = 0:1,
samples = 20,
samples y = 20,
surf,
fill=SpringGreen,
faceted color = teal] {-(x-0.5)^2/0.36 + -(y-0.5)^2/0.25 +1};
\draw[->, very thick] (0,0,0.5)node[below]{0} -- (0.93,0,0.5) node[below]{1} -- (1,0,0.5) node[above] {$\boldsymbol{u}_1$};
\end{axis}
\node[below, xshift = 2.75 cm, yshift = -0.2cm]{(a)};



\begin{axis}[hide axis, xmin=-0.2, xmax = 1.2, ymin = -0.2, ymax = 1.2, xshift = 5.5 cm,  scale =0.7]

\draw[->, thick] (0,0) -- (0,1)node[left]{$\boldsymbol{u}_2$};
\draw[teal, thick] (0.5,0.5) ellipse (0.04*0.6 and 0.04*0.5);
\draw[teal, thick] (0.5,0.5) ellipse (0.16*0.6 and 0.16*0.5);
\draw[teal, thick] (0.5,0.5) ellipse (0.36*0.6 and 0.36*0.5);
\draw[teal, thick] (0.5,0.5) ellipse (0.64*0.6 and 0.64*0.5);
\draw[teal, thick] (0.5,0.5) ellipse (0.6 and 0.5);
\draw[BrickRed, thick] (0.5-0.04*0.6,0.5-0.04*0.5) rectangle (0.5+0.04*0.6, 0.5+0.04*0.5);
\draw[BrickRed, thick] (0.5-0.16*0.6,0.5-0.16*0.5) rectangle (0.5+0.16*0.6, 0.5+0.16*0.5);
\draw[BrickRed, thick] (0.5-0.36*0.6,0.5-0.36*0.5) rectangle (0.5+0.36*0.6, 0.5+0.36*0.5);
\draw[BrickRed, thick] (0.5-0.64*0.6,0.5-0.64*0.5) rectangle (0.5+0.64*0.6, 0.5+0.64*0.5);
\draw[BrickRed, thick] (0.5-0.5,0.5-0.5) rectangle (0.5+0.5, 0.5+0.5);
\draw[->, thick] (0,0)node[below]{0} --  (1,0) node[below] {$\boldsymbol{u}_1$};

\end{axis}
\node[below, xshift = 7.95 cm, yshift = -0.2cm]{(b)};

\begin{axis}[ xmin=-0.1, xmax = 1, zmin=0, zmax=1, hide axis, xshift = 10.5 cm, scale =0.8]
\draw[->, very thick] (0,0,0.5) -- (0,1,0.5)node[left]{$\boldsymbol{u}_2$};
\addplot3 [
domain=0:1,
domain y = 0:1,
samples = 20,
samples y = 20,
surf,
fill=Maroon!30,
faceted color = BrickRed] {-max(abs(x-0.5)/0.6, abs(y-0.5)/0.5)^2 +1};
\draw[->, very thick] (0,0,0.5)node[below]{0} -- (0.93,0,0.5) node[below]{1} -- (1,0,0.5) node[above] {$\boldsymbol{u}_1$};
\end{axis}
\node[below, xshift = 13.5 cm, yshift = -0.2cm]{(c)};
    \end{tikzpicture}
    \caption{(a) Function $\psi$ from Example~\ref{ex:paraboloid}. (b) In green, contour lines of $\psi$, the interior of which are the slices $\{\psi>z_i\}$ to some $z_i\in\mathbb{R}$, $i\in[5]$. In red, contour of the sets $I_{\{\psi>z_i\}\cap[0,1]^2}$. (c) Function $\psi^*$ derived from those contour lines. We have $\int \psi \,\mathrm{d} \uC_2 = \int \psi^*\,\mathrm{d} \uC_2 = 2/3$.}
    \label{fig:exintegr}
\end{figure}
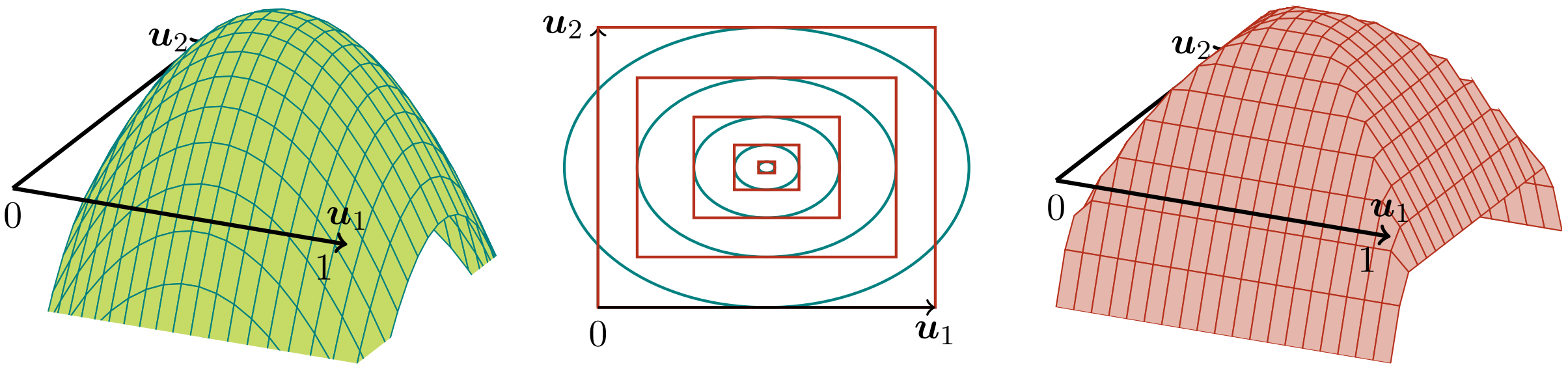

\begin{example}
\label{ex:paraboloid}
Consider the function $\psi:[0,1]^2\mapsto \mathbb{R}$ written 
\begin{equation*}
\psi(u_1,u_2) = 1 -\frac{(u_1-0.5)^2}{0.36} - \frac{(u_2-0.5)^2}{0.25}, \quad \text{for }u_1,u_2\in\mathbb{R}. 
\end{equation*}
The curve described by $\psi$ is a paraboloid centered at $(0.5,0.5)$; see Figure~\hyperref[fig:exintegr]{\ref{fig:exintegr}(a)} for an illustration. This function is concave, and hence quasiconcave on $[0,1]^2$, so Proposition~\ref{prop:integr-gravitational} applies. Its superlevel sets $\{\mathbf{u}\in[0,1]^2: \psi(\mathbf{u})> z\}$, for every $z<1$, are ellipses clipped to $[0,1]^2$ whose axial diameters are respectively $1.2\sqrt{1-z}$ and $\sqrt{1-z}$; see Figure~\hyperref[fig:exintegr]{\ref{fig:exintegr}(b)}.  
Therefore, 
\begin{equation*}
    \ukappa_1(z)-\lkappa_1(z) = \min\left(1.2\sqrt{1-z},1\right) \quad \text{and}\quad  \ukappa_2(z)-\lkappa_2(z) = \min\left(\sqrt{1-z},1\right),\quad \text{for every }z<1,  
\end{equation*}
so $\min_{i\in\{1,2\}}(\ukappa_i(z) - \lkappa_i(z)) = \min(\sqrt{1-z},1)$. We have $\inf \psi = -25/36$ on $[0,1]^2$. By Proposition~\ref{prop:integr-gravitational}, 
\begin{equation*}
\uJ_2(\psi) = -\frac{25}{36} + \int_{-25/36}^1 \min\left(\sqrt{1-z},1\right)\mathrm{d}z = -\frac{25}{36} + \frac{25}{36} + \int_0^1 \sqrt{1-z}\mathrm{d}z = \frac{2}{3}. 
\end{equation*}


\end{example}

\section{Connections to problems in various fields}
\label{sect:app}
The upper and lower copula capacities have important roles in several fields, and we list seven of them in this section. 

\subsection{Fractional matching in multipartite hypergraphs}
\label{subsect:app-matching}

Consider the following matching problem. 
There are $n$ groups each comprising $m$ individuals. The objective is to create as many teams as possible, where teams comprise one individual from each group. However, some teams are forbidden (for instance, because some individuals dislike each other.) If each individual can only be on one team, then the matching is said to be \textit{integral}; if individuals can divide their time between multiple teams, then the matching is said to be \textit{fractional}. The integral problem is exposed in \citet[Section~2, Problem~7 and Section~4, Problem~17]{K72} for two and three groups respectively.

Usually, such matching problems are represented by a multipartite hypergraph where edges connect admissible teammates, but we can represent it differently to fit our framework:
We partition $[0,1]$ into intervals $I_1,\ldots,I_m$ each of length $1/m$. Let $E\subseteq [m]^n$ be the set of admissible teams; these are the edges of the multipartite hypergraph. We define 
\begin{equation}
    A = \bigcup_{(j_1,\ldots,j_n)\in E} I_{j_1}\times \ldots \times I_{j_n}, \label{eq:A-matching} 
\end{equation}
meaning each cell from the grid-like partition of $[0,1]^n$ belonging to $A$ represents an admissible team. 

The upper copula capacity then provides the solution to the fractional matching problem: for $A$ as in \eqref{eq:A-matching}, the value $m\uC_n(A)$ is the maximum number of admissible teams that can be created. Note that, indeed, in a fractional matching, every group has the same total unused time, which can thus be coupled independently to complete the copula. 
Since every integral matching is a fractional matching,  $m\uC_n(A)$ is also an upper bound on the maximum number of teams for the integral matching problem. 
For dimension $n=2$, the fractional- and integral-matching problems have the same optimal value, and in fact \eqref{eq:hall-deficiency} is a continuous analogue to \cite{H35}'s condition for perfect matching. 
For dimensions $n\geq 3$, the existence of a perfect fractional matching may not guarantee the existence of a perfect integral, as we illustrate in the following example.  
\begin{example}
Consider three groups of two individuals, 0 and 1, each. Let the admissible teams be $\{000,011,101,110\}$. By defining $A$ as in \eqref{eq:A-matching}, we have $\uC_3(A) = 1$; thus, $2\uC_3(A) = 2$ teams is the maximal number of teams for a fractional matching.
It is feasible notably by letting each individual dedicate their time equally to each of the two teams they feature in. The maximal number of teams for an integral matching is 1, however, as any two admissible teams share an individual. 
\end{example}
%
%
%
The results from Section~\ref{sect:geometry} can prove useful to compute $\uC_n(A)$ and solve fractional matching in higher dimensions. Finally, let us remark that integral perfect matching in dimensions $n\geq 3$ is NP-complete \citep{K72}, so a computable upper bound given by its fractional matching analogue is of practical value.







\subsection{Existence of probability measures with given marginals}
\label{subsect:app-strassen}

As mentioned previously, to verify whether $\uC_n(A) = 1$ for a closed set $A\in\mathcal{B}_n$ is a special case of the marginal problem for the case where marginals are uniform. We discuss here that the upper copula capacity is also relevant to marginal problems with any marginals. 

Let $\mu_1,\ldots,\mu_n$ be prescribed marginal distributions, and $A^{\dagger}\subseteq\mathbb{R}^n$ be closed. 
The marginal problem is interested in whether there exist a distribution $\pi$ whose support is contained in $A^{\dagger}$ and with marginals $\mu_1, \ldots, \mu_n$. 
The marginal problem has previously been studied by \cite{S65}, \cite{E78}, \cite{GR84}, \cite{K84},  \cite{S93} and \cite{RR95} among others. If $\Pi$ is the set of all measures $\pi$ with the prescribed marginals, then the marginal problem becomes whether $\sup_{\pi\in\Pi} \pi(A^{\dagger}) = 1$. 

For $i\in[n]$, let $F_i^{-1}$ be the left quantile function of $\mu_i$. Then, by Sklar's theorem, 
\begin{equation}
\uC_n\left(\left\{\mathbf{u}\in[0,1]^n : \left(F_1^{-1}(u_1),\ldots,F_n^{-1}(u_n)\right)\in A^{\dagger} \right\}\right) = \sup_{\pi\in\Pi} \pi\left(A^{\dagger}\right).
\label{eq:marginal-problem}
\end{equation}
Thus, results for the upper copula capacity also apply to marginal problems with general marginals through the relation in \eqref{eq:marginal-problem}. As previously indicated, Proposition~\ref{prop:dim2} follows from \cite{K84}; other results from Section~\ref{sect:geometry} provide geometric criteria for the existence of such distributions, especially in higher dimensions, for which the characterization in \eqref{eq:strip-covers}--\eqref{eq:hall-deficiency} does not extend.

 \subsection{Joint mixability}
 \label{subsect:app-jm}

A specific kind of marginal problem is when the set $A^{\dagger}$ is such that all of its elements have constant sum. In other words, one wants to verify whether 
\begin{equation}
   \uC_n \left(\left\{\mathbf{u}\in [0,1]^n: \sum_{i=1}^n F^{-1}_i(u_i)= r\right\} \right) = 1 \label{eq:jointly-mixable}
\end{equation}
for some $r\in\mathbb{R}$. If \eqref{eq:jointly-mixable} holds, then $\mu_1,\ldots,\mu_n$ are said to be jointly mixable; and if all these distributions are the same as $\mu$, then $\mu$ is $n$-completely mixable. 
For the determination of whether distributions are jointly mixable, see  \cite{WW16} and \cite{PRWW19}. Joint mixability has application in various fields including combinatorial number theory \citep{LS26,KL26}, 
optimal transport (see Section~\ref{subsect:app-OT}), labour market \citep{BTZ25}, 
and
game theory (see Section~\ref{subsect:app-resources}).

 \subsection{Multiple hypothesis testing} 
 \label{subsect:app-hypothesis}
  Let $p_1,\dots,p_n$ be $n$ p-values testing a global null hypothesis $H_0$.   
  A test $A$ on the $n$ p-values is represented by a decreasing set  $A\in\mathcal{B}_n$,  
 such that $H_0$ is rejected if $(p_1,\dots,p_n)\in A$.
 For each  test $A$, the value $\uC_n(A)$ represents the supremum of its type-I error of the test   for dependent p-values.  Validity of the test  $A$ at level $\alpha$ is the constraint $\uC_n(A)\le \alpha$.  For instance, the test based on the Bonferroni correction at level $\alpha \in (0,1)$ chooses 
 $$
 A_{\rm Bonf}(\alpha)=\left\{\mathbf{u} \in[0,1]^n:  \min (\mathbf{u}) \le \alpha/n \right\},
 $$ 
 and one can verify $\uC_n(A_{\rm Bonf}(\alpha)) = \alpha$.
 Therefore, the capacity  $\uC_n$ is the key quantity to guarantee the validity of a test based on arbitrarily dependent p-values. The problem of testing with arbitrarily dependent p-values is studied in detail by \cite{VW20} and \cite{VWW22}.
Other settings in multiple hypothesis testing, such as controlling the probability of false discovery proportions, 
can also be expressed using $\uC_n$, although in  more complicated forms.

    \subsection{Robust risk aggregation and reliability theory}
    \label{subsect:app-robust}
    Let $X_1,\dots,X_n$ be random losses with specified respective distribution functions $F_1,\dots,F_n$ but unknown dependence structure. A classic question in robust risk aggregation is to determine 
     the worst-case excess probability 
    $$
      \sup_{X_i\sim F_i,~i\in [n]} \p\left(g(X_1,\dots,X_n)> x\right),$$
      where $g:\R^n\to \R$ represents a measurable aggregation function, and $x$ is a constant,
      or to determine 
      the worst-case Value-at-Risk  (VaR) 
    $$
  \sup_{X_i\sim F_i,~i\in [n]} \VaR_\alpha\left(g(X_1,\dots,X_n)\right) ,
    $$ 
    where  $\alpha\in (0,1)$, and $\VaR_\alpha$  represents the left $\alpha$-quantile of a random variable.
    The value of  the worst-case excess probability  is equal to 
    $\uC_n(A)$ where 
    $$
    A=\left\{\mathbf{u}\in [0,1]^n: g\left(F^{-1}_1(u_1),\dots,F^{-1}_n(u_n)\right)>x\right\}.
    $$
    The value of the worst-case VaR can be obtained by inverting the worst-case excess probability over $x\in \R$.
    Computing the worst-case VaR or excess probability is a difficult question by itself, and it has been actively studied in risk management; see \cite{WPY13}, \cite{EPR13, EWW15}, and the references therein, where most results are obtained for the linear aggregation function 
    $g:\mathbf{y}\mapsto  \mathbf{w} \cdot \mathbf{y}$ for some $\mathbf{w} \in \R^n$, and \cite{BLLW24} for more recent progress on this problem.
    A similar problem appears in reliability theory where the quantity of interest is the probability of communication error; see \cite{BJ20}.

      \subsection{Multi-marginal optimal transport}
      \label{subsect:app-OT}
Multi-marginal optimal transport problems are widely studied in economic theory and mathematical optimization; see \cite{S15} for an introduction. 
In this context, 
the objective is typically to minimize the expected cost
  $$
  \int_{\R^n} c (\mathbf x) \pi (\d  \mathbf x) 
  $$
  over $\pi\in \Gamma(F_1,\dots,F_n)$, where $\Gamma(F_1,\dots,F_n)$ is the set of distributions on $\R^n$ with marginals $F_1,\dots,F_n$, and $c:\R^n\to \R$ is a cost function. The multi-marginal optimal transport problem is generally difficult to solve. 
  Several special cases of cost functions are closely related to  $\uC_n$ and $\lC_n$, and we give two examples below.
  The first is an indicator loss function $c = \id_{A^{\dagger}} $ for some Borel $A^{\dagger}\subseteq\mathbb{R}^n$,
  and the infimum of the expected cost is $\lC_n(A) =1- \uC_n(A^c)$ where $A$ is defined as in \eqref{eq:marginal-problem}. For instance, if $n=2$, we may be interested in the  $0$-$1$ cost $\id_{A^{\dagger}}$ for $A^{\dagger}=\{(x,y)\in \mathbb{R}^2: x\ne y\}$ or 
   $A^{\dagger}=\{(x,y)\in \mathbb{R}^2:| x- y|>\epsilon\}$ for some $\epsilon>0$.
  The second is the  repulsive harmonic cost function  
  $c:(x_1,\dots,x_n)\mapsto -\sum_{i,j\in [n]} (x_i-x_j)^2$, 
  and it is connected to the problem of joint mixability from Section~\ref{subsect:app-jm}. 
  Indeed, $-\sum_{i,j\in [n]} (x_i-x_j)^2 = 2(\sum_{i} x_i)^2 - 2n\sum_{i=1}^n x_i^2 $, so by Jensen's inequality a joint mix guarantees an optimal transport, provided the distributions have finite second moment.

  \subsection{Resource allocation games and all-pay contests}
  \label{subsect:app-resources}
  In a limited resources allocation game, each agent allocates resources to  $n$ locations  to maximize the expected payoff from the game. 
The total resource available is a constant and the allocations can be random.
  A representative game of this kind is the Colonel Blotto game (e.g., \cite{R06}), which is a two-agent zero-sum game, and a special type of all-pay contest; see \cite{S09}.  
A standard method to analyze a Nash equilibrium strategy (for each agent) in this type of games is first calculating a potential equilibrium distribution $F_i$ of the allocated resources at each location $i\in [n]$ by a Lagrange multiplier method, and then constructing a joint distribution with marginals $F_1,\dots,F_n$ such that the sum of the allocated random resources does not exceed   the constant $s$  which is the amount of the total resource of the agent (typically the sum is equal to $s$). Whether the second step is viable is a non-trivial question, which can be formulated using $\uC_n$, as the existence of the joint distribution in the second step is equivalent to
    $\uC_n(A_{s})=1$ where 
 $$
    A_{s}=\left\{\mathbf{u}\in [0,1]^n: \sum_{i=1}^n F^{-1}_i(u_i)\le s\right\}.
 $$
 In the case where each resources are allocated uniformly on some locations, then $F_i^{-1}(u)=\lambda_i u$ for all $u \in(0,1)$, and some $\lambda_i>0$, $i\in[n]$ and the set above is precisely the one from Theorem~\ref{th:weighted-uniform-sums}. 

 \section{Conclusion}
 \label{sect:conclusion-copula}

A main objective of this paper, in introducing the upper copula capacity, was to provide a unified framework for many problems related to dependence uncertainty. As we have discussed in Section~\ref{sect:app}, such problems arise among various disciplines,  including  optimal transport, risk management, statistics, combinatorics, and decision theory. The unified framework given by the upper copula capacity thus offers a bridge for results stemming from one strand of research to another.  Theorem~\ref{th:weighted-uniform-sums}, for example, was inspired by a result in the theory of joint mixability, while Proposition~\ref{prop:dim2} comes from duality theory. 
Section~\ref{sect:INT} concerns problems in decision theory via the Choquet expected utility theory of \cite{S89}.
In this vein, Examples~\ref{ex:octahedron} and~\ref{ex:inscribed-octahedron-exact}, featuring the octahedron, and the geometric observation extracted from juxtaposing Theorem~\ref{th:convex} and Proposition~\ref{prop:dim2},  related to \cite{B51}'s plank conjecture, are illustrations of how the results above can combine to provide various insights on these problems.







\appendix

\section{Proofs}

\subsection{Proof of Proposition~\ref{prop:capacity}}

  It is straightforward to check that $\uC_n$ is a normalized capacity.
For item~\ref{item:property-2}: For any $C\in \mathbb C_n$ and $A,B\in\mathcal{B}_n$, we have 
 $ C(A\cup B)  \le C(A) + C(B)  
  $; taking suprema gives subadditivity of $\uC_n$.  For item~\ref{item:property-3}: from the increasingness of the sets, the limit can be replaced by a supremum; suprema are commutative, and every copula measure is continuous from below. In fact, all upper probabilities are capacities satisfying \ref{item:property-1}, \ref{item:property-2} and \ref{item:property-3}; see e.g.~\cite{HS73}. Because $\mathbb{C}_n$ is weakly compact, $\uC_n$ satisfies \ref{item:property-4} by \citet[Lemma 2.3]{HS73}. 
  \hfill$\square$


\subsection{Proof of Proposition~\ref{prop:max-is-attained}}

  Take a closed set $A$ and  a sequence $(C_m)_{m\in \N}$ in $\mathbb C_n$ such that $C_m(A)\to \uC_n(A)$ as $m\to \infty$. Since any element of $\mathbb C_n$ is supported in the compact set $[0,1]^n$,
  the sequence we took has a weakly converging subsequence, which  converges to $C^*$. We have $C^*\in \mathbb C_n$ because the weak limit of copulas is again a copula.
  It follows that $C^*(A)=\uC_n(A)$ since $A$ is closed, and thus the maximum is attained. \hfill $\square$

\subsection{Proof of Proposition~\ref{prop:lines}}

By Proposition~\ref{prop:max-is-attained}, the set $\mathcal{M} = \{C\in\mathbb{C}_n : C(A) = \uC_n(A)\}$ is nonempty since $A$ is closed. The set $\mathcal{M}$ is also weakly closed: Consider any sequence $\{C_k : k\in\mathbb{N}\}\subseteq \mathcal{M}$ weakly converging to $C\in\mathbb{C}_n$. The Portmanteau theorem yields $C(A) \geq \limsup_{k\to\infty} C_k(A) = \uC_n(A)$, so $C\in\mathcal{M}$ as well. Because $\mathbb{C}_n$ is weakly compact, weak closedness of $M$ means it is also weakly compact. 

For any nonempty compact convex set, the Krein-Milman theorem guarantees the existence of extreme points in $\mathcal{M}$. Let $C^*$ be an extreme point in $\mathcal{M}$; we show that it is also  an extreme point in $\mathbb{C}_n$. Suppose it is not: then, there are $C_1,C_2\in\mathbb{C}_n$, different from $C^*$, such that $C^*=tC_1 + (1-t)C_2$ for some $t\in(0,1)$. But because $C^*(A) = \uC_n(A)$, this means that $C_1(A) = C_2(A) = \uC_n(A)$ because $\uC_n(A)$ is the largest value attained by a copula. Hence, $C_1,C_2\in\mathcal{M}$, and this would contradict the extremality of $C^*$ in $\mathcal{M}$. 
\citet[Theorem 2.1]{S86} indicates that $C^*\perp \ell_n$ because it is an extreme point in $\mathbb{C}_n$. \hfill$\square$

\subsection{Proof of Proposition~\ref{prop:intervals}}

By monotonicity, 
\begin{equation*}
    \uC_{k+m}(A\times B) \leq \uC_{k+m}\left(A\times[0,1]^m\right) = \uC_k(A),  
\end{equation*}
where the last equality holds because the set of $k$-marginals of higher dimensional copulas is the same as $\mathbb{C}_k$. 
Similarly, $\uC_{k+m}(A\times B) \leq \uC_m(B)$, so combining both gives $\uC_{k+m}(A\times B) \leq \min(\uC_k(A),\uC_m(B))$. 
For the converse inequality, choose any $C^{\prime}\in\mathbb{C}_k$ and $C^{*}\in\mathbb{C}_m$. Let $q^{\prime} = C^{\prime}(A)$ and $q^* = C^{*}(B)$. By symmetry, we can suppose $q^{\prime}\leq q^*$. If $q^{\prime} = 1$ or $q^{*} = 0$, the product measure $C^{\prime}\otimes C^*$ already gives the desired value. Otherwise, define the set function $\mu$ given by $\mu(E) = C^{*}(E\cap B)/q^*$ for every $E\in\mathcal{B}_m$. 
It is straightforward to see that $\mu$ is a propbability measure concentrated on $B$, and $C^* - q^{\prime}\mu$ is nonnegative because $q^{\prime} \leq q^*$. We write $C_{|A}^{\prime}(E) = C^{\prime}(A\cap E)$. 
Define the measure $C^{\dagger}$ by 
\begin{equation*}
    C^{\dagger} = C^{\prime}_{|A}\otimes \mu + C^{\prime}_{|A^{\rm c}}\otimes \left(\frac{C^* - q^{\prime}\mu}{1-q^{\prime}} \right). 
\end{equation*}
The measure $C^{\dagger}$ is a copula: its first $k$ coordinates have $k$-marginal $C^{\prime}_{|A} + C^{\prime}_{|A^{\rm c}} = C^{\prime}$, its last $m$ coordinates have $m$-marginal $q^{\prime}\mu + (C^{*} - q^{\prime}\mu) = C^*$, 
and therefore all coordinates have uniform marginals. 
We have $C^{\dagger}(A\times B) = q^{\prime} = \min(C^{\prime}(A), C^*(B))$. Since this holds for any choice of $C^{\prime}$ and $C^*$, we have $\min(\uC_k(A), \uC_m(B))\leq \uC_{m+k}(A\times B)$.  

Combining the two inequalities gives the first part of the proposition; 
induction using $m=k=1$ directly yields the second part of the proposition. 
\hfill$\square$

\subsection{Proof of Theorem~\ref{th:rearrangements}}

Write $B = \{\mathbf{u}\in[0,1]^n : (\varphi_1(u_1),\ldots,\varphi_n(u_n))\in A\}$. Choose any $C\in\mathbb{C}_n$. For any random vector $\mathbf{U}\sim C$, the random vector $(\varphi_1(U_1),\ldots,\varphi_n(U_n))$ also has uniform marginals by definition of rearrangements of the identity. Thus, 
\begin{equation*}
    C(B) = \mathbb{P}\left( (\varphi_1(U_1),\ldots,\varphi_n(U_n)) \in A\right) \leq \uC_n(A), 
\end{equation*}
and since this holds for any $C\in\mathbb{C}_n$, we have $\uC_n(B) \leq \uC_n(A)$. 

For the converse inequality, choose any $C^*\in\mathbb{C}_n$ and consider a random vector $\mathbf{V}\sim C^*$. We search for a copula $C^{\dagger}$ such that $\varphi_i(U_i) = V_i$ almost surely for some $\mathbf{U}\sim C^{\dagger}$. For every $i\in[n]$, define the kernel $\kappa_i$ by letting $\kappa_i(x, E) = \mathbb{P}(W\in E\mid \varphi_i(W) = x)$ where $W$ has a uniform distribution on $(0,1)$, for every $x\in[0,1]$ and $E\in\mathcal{B}$. Then, for every $i\in[n]$, let $(U_i,\mathbf{V})$ have the joint distribution given by
\begin{equation*}
    \mathbb{P}(U_i\in E, \mathbf{V}\in E^{\prime}) = \int_{E^{\prime}} \kappa_i(x_i, E)C^*(\mathrm{d}x), \quad \text{for every }E\in\mathcal{B},\,E^{\prime}\in\mathcal{B}_n. 
\end{equation*}
By letting the $U_i$ be conditionally independent given $\mathbf{V}$, we obtain a random vector $\mathbf{U}$ such that $\varphi_i(U_i) = V_i$ almost surely, and the marginal distributions of $\mathbf{U}$ are uniform because $\varphi_i$, $i\in[n]$, are rearrangements of the identity. Take $C^{\dagger}$ as the distribution of $\mathbf{U}$. As a result, 
\begin{equation*}
    C^*(A) = \mathbb{P}(\mathbf{V}\in A) = \mathbb{P}((\varphi_i(U_i):i\in[n]) \in A) = \mathbb{P}(\mathbf{U}\in B) = C^{\dagger}(B) \leq \uC_n(B).
\end{equation*}
Since this holds for any $C^*\in\mathbb{C}_n$, we have $\uC_n(A) \leq \uC_n(B)$. We conclude that $\uC_n(A) = \uC_n(B)$. \hfill$\square$

\subsection{Proof of Corollary~\ref{prop:translations}}

For any $\mathbf{x}\in\mathbb{R}^n$, the equivalence $\mathbf{u}\in\tau_{\boldsymbol{x}}(A) \iff \langle \mathbf{u} - \mathbf{x}\rangle \in A$ holds for almost every $\mathbf{u}\in [0,1)^n$, for every copula. 
The function $u_i\mapsto \langle u_i - x_i\rangle$ belongs to $\Phi_U$. The claim thus follows by Theorem~\ref{th:rearrangements}. 
\qed

\subsection{Proof of Proposition~\ref{lem:closure}}

Note that $\overline{A^{\circ}}\subseteq \overline{A}$. The condition that $A$ is strictly star-shaped means that we have $t\overline{A} + (1-t)\mathbf{a}\subseteq A^{\circ}$ for every $t\in(0,1)$, for some $\mathbf{a}\in A^{\circ}$. By letting $t\uparrow 1$, we have $\overline{A} = \overline{A^{\circ}}$. By Proposition~\ref{prop:max-is-attained}, there is a copula $C^* \in\mathbb{C}_n$ satisfying $C^*(\overline{A}) = \uC_n(\overline{A})$. Consider a random vector $\mathbf{U}\sim C^*$. For any $t\in(0,1)$, the random vector $(1-t)\mathbf{a} + t\mathbf{U}$ has marginal distributions uniform on $J_i = ((1-t)a_i, (1-t)a_i+t)$ for each $i\in[n]$. Consider the random vector $\mathbf{V}$ with mutually independent components, each following a uniform distribution on $[0,1]\backslash J_i$ for every $i\in[n]$. Then, the mixture in respective proportion $t$ and $(1-t)$ of the distributions of $t\mathbf{U}+(1-t)\mathbf{a}$ and $\mathbf{V}$ is a copula, and let $C_t$ denote this copula for every $t\in(0,1)$. We thus have 
\begin{equation}
    C_t(A^{\circ}) = t\mathbb{P}(t\mathbf{U}+(1-t)\mathbf{a} \in A^{\circ}) + (1-t)\mathbb{P}(\mathbf{V}\in A^{\circ}). 
    \label{eq:closure-1}
\end{equation}
Moreover, because $t\overline{A} + (1-t)\mathbf{a}\subseteq A^{\circ}$, 
\begin{equation}
    C^*(\overline{A}) \leq \mathbb{P}(t\mathbf{U}+(1-t)\mathbf{a} \in A^{\circ})\label{eq:closure-2}.
\end{equation}
Combining \eqref{eq:closure-1} and \eqref{eq:closure-2} gives 
\begin{equation*}
    t\uC_n(\overline{A}) \leq tC^*(\overline{A}) \leq  C_t(A^{\circ}) \leq \uC_n(A^{\circ}) \quad \text{for every }t\in(0,1).
\end{equation*}
Letting $t\uparrow 1$ yields $\uC_n(\overline{A}) \leq \uC_n(A^{\circ})$. Finally, let us remark that $A^{\circ}\subseteq A \subseteq \overline{A}=\overline{A^{\circ}}$; the desired result follows by monotonicity. \hfill$\square$

\subsection{Proof of Proposition~\ref{prop:dim2}}
\label{proof:dim2}

The identity in \eqref{eq:strip-covers} is obtained directly by combining Theorem~2.16 and Proposition~3.3 of \cite{K84}; note that $\uC_2$ coincides with his functional $S$ for indicator functions (which are upper-semicontinuous for closed sets) under uniform marginals. We argue below the equivalence between \eqref{eq:strip-covers} and \eqref{eq:hall-deficiency}. 

For every closed set $K\subseteq [0,1]$, we can choose sets $S^* = [0,1]\backslash K$ and $T^* = \pi_A(K)$, and they satisfy $A\subseteq p_1^{-1}(S^*) \cup p_2^{-1}(T^*)$ by definition of $\pi_A(K)$. We have $\ell(S^*) + \ell(T^*) = 1- \ell(K) + \ell(\pi_A(K))$. Thus, 
\begin{equation*}
    \inf_{\substack{A\subseteq(p_1^{-1}(S)\cup p_2^{-1}(T))\\ S,T\in\mathcal{B}}}\{\ell(S)+\ell(T)\}
 \leq \inf_{\substack{K\subseteq[0,1]\\K \text{ is closed}}}\{ 1- \ell(K)+\ell(\pi_A(K))\}. 
\end{equation*}
For the converse inequality, choose any $S,T$ such that $A\subseteq p_1^{-1}(S) \cup p_2^{-1}(T)$. This means that, for any $(x,y) \in A$, if $x\not\in S$ then $y\in T$, and therefore 
\begin{equation}
    \pi_A([0,1]\backslash S)\subseteq T. \label{eq:pia-k}
\end{equation}
For any $\varepsilon>0$, there exists a closed $K^*$ satisfying $K^*\subseteq [0,1]\backslash S$ and $\ell(K^*) \geq 1- \ell(S)-\varepsilon$. By \eqref{eq:pia-k}, we have $\pi_A(K^*)\subseteq T$, and therefore 
\begin{equation}
     \inf_{\substack{K\subseteq[0,1]\\K \text{ is closed}}}\{ 1- \ell(K)+\ell(\pi_A(K))\} \leq 1-\ell(K^*) + \ell(\pi_A(K^*)) \leq \ell(S) + \varepsilon + \ell(T).
\end{equation}
Taking $\varepsilon\to 0$ yields the desired inequality. 

The last assertion of the proposition is supported by Proposition~\ref{prop:max-is-attained}, ensuring the existence of a copula supported on $A$ if $\uC_2(A) = 1$ and $A$ is closed. \hfill$\square$

\subsection{Proof of Theorem~\ref{th:convex}}

The second equality, $\uC_2(I_A) = \min_{i\in\{1,2\}}(b_i(A) - a_i(A))$ is given by Proposition~\ref{prop:intervals}. We need to argue $\uC_2(A) = \uC_2(I_A)$. We have $\uC_n(A)\leq \uC_n(I_A)$ directly by monotonicity, and this also makes equality straightforward if either projection of $I_A$ is a null set. Thus, we can assume  $\min_{i\in\{1,2\}}(b_i(A) - a_i(A))>0$ for the remainder, and we want to show that $\uC_2(A)\geq \uC_2(I_A)$.

First, we assume that $A$ is closed. Define the set
\begin{equation*}
    B = \left\{ \left(\frac{u_1 - a_1(A)}{b_1(A) - a_1(A)}, \frac{u_2-a_2(A)}{b_2(A) - a_2(A)} \right):(u_1,u_2)\in A \right\}. 
\end{equation*}
It follows that $B$ is closed, convex and touches all four edges of $[0,1]^2$. By Lemma~\ref{lem:copula-construction} (Appendix~\ref{sect:construction-copula}), there exists a random vector $(V_1,V_2)$ with uniform marginals and whose support is contained in $B$. Let $X_i = a_i(A) + (b_i(A) - a_i(A)) V_i$ for each $i\in\{1,2\}$. By the definition of $B$, we have $(X_1,X_2)\in A$ almost surely; also, $X_i$ follows a uniform distribution on $(a_i(A),b_i(A))$ since $V_i$ follows a uniform distribution on $(0,1)$, for each $i\in\{1,2\}$.
Consider any copula $C\in\mathbb{C}_2$, and let $(U_1,U_2)\sim C$ be independent of $(X_1,X_2)$. Define $(W_1,W_2)$ by letting
\begin{equation}
    W_i = \begin{cases}
    X_i,\quad \text{if }U_i\in[a_i(A),b_i(A)],\\
        U_i,  \quad \text{if }U_i\not\in[a_i(A),b_i(A)]
    \end{cases}, \quad \text{for } i\in\{1,2\}.
\end{equation}
Because $X_i$ is uniformly distributed on $(a_i(A),b_i(A))$, then $(W_1,W_2)$ has  marginal distributions uniform on $(0,1)$.  Let $C^{\dagger}$ be the distribution of $(W_1,W_2)$; note that $C^{\dagger}(I_A) = C(I_A)$. Also, $C^{\dagger}(I_A) = C^{\dagger}(A)$ because $(X_1,X_2)$'s support is contained in $A$. 
Thus, $C(I_A) = C^{\dagger}(A) \leq \uC_n(A)$. Because this holds for all $C\in\mathbb{C}_2$, we have $\uC_2(I_A) \leq \uC_2(A)$, and hence $\uC_2(I_A) = \uC_2(A)$.

If $A$ is not closed and has non-empty interior, then it is strictly star-shaped due to its convexity. We can apply the argument above for $\overline{A}$, and $\uC_2(A) = \uC_2(\overline{A})$ by Proposition~\ref{lem:closure}. If $A$ is not closed and has empty interior, then it is a segment and every copula will yield the same value whether its endpoints are included or not.   \hfill$\square$

\subsection{Proof of Proposition~\ref{th:decreasing-sets}}
  Let $\mathcal U$ be the set of all random variables  $U$ with a uniform distribution on $(0,1)$, 
  and $\mathcal V$ be the set of all random variables  $V$ with $\p(V\le t)\le t$ for all $t\in [0,1]$.  In other words, $\mathcal{V}$ is the set of random variables whose distributions stochastically dominate the uniform distribution on $(0,1)$. For $\lambda = 0$, the set $\lambda A$ is either empty or a singleton, and the claim is then straightforward; assume below that $\lambda\in(0,1]$
The following lemma is useful in proving the result, and it follows directly from Theorem 6.B.1 of \cite{SS07}.
\begin{lemma}  \label{lem:dominated} For any decreasing set $A\in\mathcal{B}_n$ and any $\mathbf V\in \mathcal V^n$, 
there exists $\mathbf U\in \mathcal U^n$ such that 
$\p(\mathbf U\in A) \ge  \p(\mathbf V\in A)$.
\end{lemma} 
  We will prove two statements, which will be sufficient for the proposition. 
  \begin{enumerate}[label={\rm (\roman*)}, ref={\rm (\roman*)}]
 \item \label{item:lem-UV-a} For any $\mathbf U\in \mathcal U^n$ and $\lambda \in (0,1]$,  
 there exists $\mathbf V\in \mathcal U^n$ such that 
 $\p(\mathbf V\in \lambda A)\ge  \lambda  \p(\mathbf U\in   A)$.
 \item  \label{item:lem-UV-b} For any $\mathbf U\in \mathcal U^n$ and $\lambda \in (0,1]$,   
 there exists $\mathbf V\in \mathcal U^n$ such that 
 $  \p(\mathbf V\in   A)\ge \p(\mathbf U\in \lambda A) /\lambda $.
 \end{enumerate}
 Statement \ref{item:lem-UV-a} implies 
$\uC_n(\lambda A) \ge \lambda \uC_n(  A)$ 
and  statement \ref{item:lem-UV-b} implies $\uC_n(\lambda A) \le \lambda \uC_n(  A)$.
 %
   
   We first show \ref{item:lem-UV-a}. 
Take any  event $A$ with probability $\lambda$ independent of $\mathbf U$ (the existence of $A$ is safe to assume, since we are interested in statements on distributions). 
Define $\mathbf V^*  $ by $$\mathbf V^*= 
\lambda \mathbf U \times  \id_A+ \mathbf 1 \times \id_{A^c}.$$  
It is straightforward to check $\mathbf V^*\in \mathcal V^n$. Hence,   
by  Lemma \ref{lem:dominated}, there exists $\mathbf V\in \mathcal U^n$ such that 
$$\p( \mathbf V   \in  \lambda  A  )  \ge \p( \mathbf V^*   \in  \lambda  A  ) \ge\p(A) \p\left (    \mathbf U\in    A  \right)= \lambda  
\p\left (    \mathbf U\in     A \right) .$$ 
Next  we show \ref{item:lem-UV-b}.   
Note that $\p(\mathbf U\in \lambda A) \le \uC_n(\lambda A )\le \uC_n([0,\lambda]^n)=\lambda$, where the last equality follows by Proposition \ref{prop:intervals}. 
Let $E$ be an event containing $\{\mathbf U\in \lambda A\}$ with probability $  \lambda  $. Let $\mathbf V^*=(V_1^*,\dots,V_n^*)$ follow the conditional distribution of $\mathbf U/\lambda$ given $E$.
Note that  for $i\in[n]$, 
$$
\p(V_i^* \le  t ) = \p\left(\left.\frac{U_i}{\lambda} \le t  \;\right|\; E\right) \le \frac{\p(U_i\le \lambda t )}{\p(E)} 
=t,
$$
and hence   $\mathbf V^*\in \mathcal V^n$. 
Moreover,  by  Lemma \ref{lem:dominated}, there exists $\mathbf V\in \mathcal U^n$ such that 
$$
\p( \mathbf V  \in A) \ge  \p( \mathbf V^* \in A) = \p(\mathbf  U \in \lambda A\mid E) =   \frac{\p(\mathbf  U \in \lambda A)}{\lambda  } .
$$  
Hence, \ref{item:lem-UV-a} and \ref{item:lem-UV-b} both hold true.\hfill$\square$

\subsection{Proof of Theorem~\ref{th:weighted-uniform-sums}}

Write $A_s = \left\{\mathbf{u}\in[0,1]^n:\sum_{i=1}^n \lambda_i u_i \leq s\right\}$. For $s = 0$, the result is immediate because $A_s = \{0\}$, and $\uC_n(\{0\}) = 0$. Thus assume that $s>0$. 

We first show that $\uC_n(A_s) \leq \min(s/L(\boldsymbol{\lambda}), 1)$. 
Choose any $C\in\mathbb{C}_n$, and consider a random vector $\mathbf{U}\sim C$. Write $p = \mathbb{P}(\mathbf{U}\in A_s)$. For each $i\in[n]$, it holds that $\{\mathbf{U}\in A_s \}\subseteq \{U_i\leq s/\lambda_i\}$, so
\begin{equation}
    p \leq \frac{s}{\lambda_i} \quad \text{for all }i\in[n].
    \label{eq:jointmix-1}
\end{equation}
Moreover, since $U_i$ follows a uniform (0,1) distribution, 
\begin{align*}
    \mathbb{E}[U_i \id_{\{\mathbf{U}\in A_s\}}] = \int_0^1 \mathbb{P}(U_i>r, U\in A_s) \mathrm{d}r \geq \int_0^{p}(p-\mathbb{P}(U_i\leq r))\mathrm{d}r = \int_0^p (p-r) \mathrm{d}r = \frac{p^2}{2}. 
\end{align*}
Hence, 
\begin{equation}
    \frac{p^2}{2} \sum_{i=1}^n \lambda_i \leq \mathbb{E}\left[\sum_{i=1}^n\lambda_iU_i\id_{\{\mathbf{U}\in A_s \}} \right] \leq ps
    \label{eq:jointmix-2}
\end{equation}
where the second inequality comes from the fact that $\sum_{i=1}^n\lambda_i U_i \leq s$ on the event $\{\mathbf{U}\in A_s\}$, and this event has probability $p$.  
From \eqref{eq:jointmix-2}, we deduce $p\leq s/(1/2\sum_{i=1}^n \lambda_i)$, which together with \eqref{eq:jointmix-1} gives
\begin{equation*}
    p \leq \frac{s}{\max( \max_{i\in[n]} \lambda_i, \frac{1}{2}\sum_{i=1}^n \lambda_i )};
\end{equation*}
this is $C(A_s) \leq s/L(\boldsymbol{\lambda})$. Also, $C(A_s) \leq 1$ evidently. Since this holds for every $C\in\mathbb{C}_n$, we have $\uC_n(A_s) \leq \min(s/L(\boldsymbol{\lambda}), 1)$. 

For the converse inequality, we first need to argue that $\uC_n(A_{L(\boldsymbol{\lambda})}) = 1$. 
If $\max_{i\in[n]} \lambda_i \leq \frac{1}{2}\sum_{i=1}^n \lambda_i$, the joint mixability theorem from \citet[Theorem 3.1]{WW16} indicates that there exists a copula $C^*$ such that for $\mathbf{U}\sim C^*$, we have $\mathbb{P}(\sum_{i=1}^n \lambda_i U_i = L(\boldsymbol{\lambda})) = 1$; thus, $C^*(A_{L(\boldsymbol{\lambda})}) = 1$ and this, in turn, indicates $\uC_n(A_{L(\boldsymbol{\lambda})}) = 1$. If $\max_{i\in[n]}\lambda_i > \frac{1}{2}\sum_{i=1}^n \lambda_i$, let $j$ be the index with $\lambda_j = \max_{i\in[n]} \lambda_i$, and thus $\lambda_j = L(\boldsymbol{\lambda})$. Define $\mathbf{V}$ as a random vector with marginal distributions uniform on $(0,1)$ with the following stochastic representation:
\begin{equation*}
    V_j = W \quad \text{and}\quad V_i = 1-W \quad \text{for all }i\neq j, i\in[n], 
\end{equation*}
where $W$ itself follows a uniform distribution on $(0,1)$. Then, 
\begin{equation*}
    \sum_{i=1}^n \lambda_i V_i = \lambda_j W + (1-W)\sum_{i\neq j} \lambda_i \leq \lambda_j, 
\end{equation*}
with the last inequality because $\sum_{i\neq j}\lambda_i < \lambda_j$. Then, if $C^{\dagger}$ is the distribution of $\mathbf{V}$, we have $C^{\dagger}(A_{L(\boldsymbol{\lambda})}) = \mathbb{P}(\mathbf{V}\in A_{L(\boldsymbol{\lambda})}) = \mathbb{P}(\sum_{i=1}^n \lambda_i V_i \leq \lambda_j) = 1$.

Hence $\uC_n(A_{L(\boldsymbol{\lambda})}) = 1$. 
Next, for any $s\geq L(\boldsymbol{\lambda})$, we have $\uC_n(A_s) = 1$ because $A_{L(\boldsymbol{\lambda})} \subseteq A_s$. For $0\leq s < L(\boldsymbol{\lambda})$, we have $(s/L(\boldsymbol{\lambda}))A_{L(\boldsymbol{\lambda})} \subseteq A_s$; thus, 
\begin{equation*}
    \uC_n(A_s) \geq \uC_n\left( \frac{s}{L(\boldsymbol{\lambda})}A_{L(\boldsymbol{\lambda})}\right) = \frac{s}{L(\boldsymbol{\lambda})} \uC_n\left(A_{L(\boldsymbol{\lambda})}\right) = \frac{s}{L(\boldsymbol{\lambda})}.
\end{equation*}
Combining both cases gives $\uC_n(A_s) \geq \min(s/L(\boldsymbol{\lambda}), 1)$. 
\hfill $\square$

\subsection{Proof of Proposition~\ref{prop:JK}}

The function $\psi$ is bounded, thus $a\leq \psi \leq b$ for some $a,b\in\mathbb{R}$. By definition of the Choquet integral, 
\begin{equation}
    \uJ_n(\psi) = a+ \int_a^b \uC_n(\psi>z) \mathrm{d}z, \quad \text{and} \quad \int\psi\mathrm{d}C = a+\int_a^bC(\psi>z)\mathrm{d}z\quad \text{ for every }C\in\mathbb{C}_n. 
    \label{eq:proofJK-choquet}
\end{equation}
For the implication \ref{item:JK-eq}$\implies$\ref{item:JK-att}: We choose a sequence $\{C_i : i\in\mathbb{N}\}\subseteq\mathbb{C}_n$ such that $\int\psi \mathrm{d}C_i \to \uK_n(\psi)$ as $i\to\infty$. Because $\mathbb{C}_n$ is weakly compact, there is a weakly converging subsequence $\{C_{i_j}:j\in\mathbb{N}\}$ converging to some $\mathbb{C}^*\in\mathbb{C}_n$. Since $\psi$ is bounded and upper semicontinuous, the Portmanteau theorem yields 
\begin{equation*}
    \uK_n(\psi) = \limsup_{j\to\infty} \int \psi \mathrm{d}C_{i_j} \leq \int \psi\mathrm{d}C^* \leq \uK_n(\psi),
 \end{equation*}
with the last inequality holding by definition of $\uK_n$. Therefore, the supremum defining $\uK_n(\psi)$ is attained by $C^*$. Then, from \eqref{eq:proofJK-choquet} and the assumption $\uJ_n(\psi)= \uK_n(\psi)$, 
\begin{equation*}
    \int_a^b\left(\uC_n(\psi>z) - C^*(\psi>z) \right)\mathrm{d}z = \uJ_n(\psi) - \int\psi\mathrm{d}C^* = \uJ_n(\psi) - \uK_n(\psi) = 0, 
\end{equation*}
The integrand in the left-hand integral must be nonnegative by definition of $\uC_n$, so this means $\uC_n(\psi>z) = C^*(\psi>z)$ for almost every $z\in(a,b)$. Because $a\leq \psi \leq b$ and by continuity from below of $\uC_n$, this identity actually holds true for all $z\in\mathbb{R}$. 

For the implication \ref{item:JK-att}$\implies$\ref{item:JK-eq}: By \eqref{eq:proofJK-choquet}, we directly have $\uJ_n(\psi) = \int\psi \mathrm{d}C^*$.  Moreover, $\int\psi \mathrm{d}C^* \leq \uK_n(\psi)$ by definition of $\uK_n$, and $\uK_n(\psi) \leq \uJ_n(\psi)$ from \eqref{eq:layers}. Hence, $\uJ_n(\psi) = \uK_n(\psi)$.   \hfill$\square$

\subsection{Proof of Proposition~\ref{prop:integr-gravitational}}

By definition of quasiconcavity, the sets $\{\mathbf{u}\in[0,1]^2 : \psi(\mathbf{u}) > z\}$ are convex and nonempty for every $z<\sup\psi$. Therefore, $\uC_2(\psi> z) = \min_{i\in\{1,2 \}} (\ukappa_i(z) - \lkappa_i(z))$ by Theorem~\ref{th:convex}. The integration formula is then derived directly from the definition of the Choquet integral; see \eqref{eq:def-choquet}. Note that, by quasiconcavity, $\inf \psi$ is finite. 
\hfill$\square$


\section{Copula construction for convex sets in dimension 2}
\label{sect:construction-copula}

\begin{lemma}
\label{lem:copula-construction}
    Any closed convex set touching the four edges of $[0,1]^2$ contains the support of a copula. 
\end{lemma}

\begin{proof}
We derive a general explicit construction for such a copula for any closed convex set $A\in\mathcal{B}_2$.

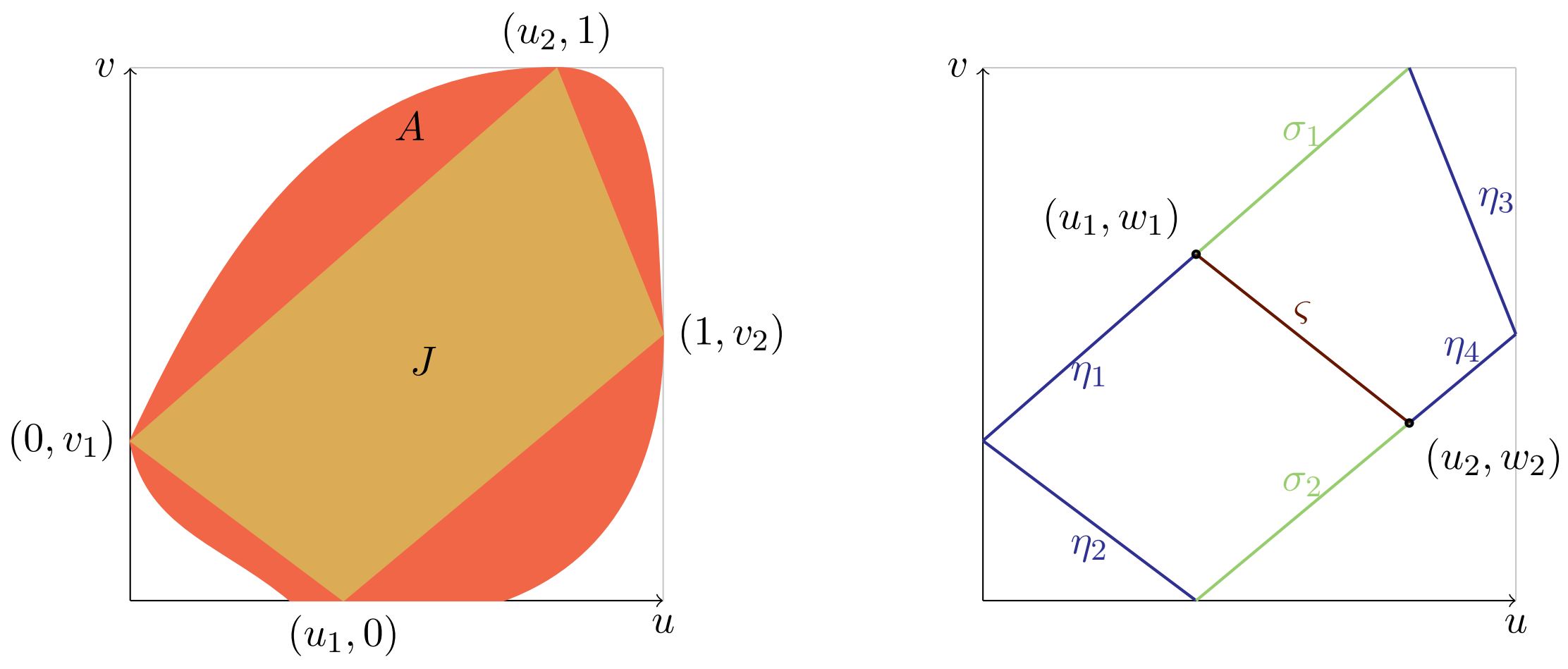
\begin{figure}[t]
\centering
    \begin{tikzpicture}[scale = 5]
          \draw[->] (0,0) -- (1,0) node[below]{$u$} ;
		\draw[->] (0,0) -- (0,1) node[left]{$v$} ;
        \draw[Black!25] (1,0) -- (1,1);
        \draw[Black!25] (0,1) -- (1,1);

            \filldraw[Red!75] (0,0.3) to [out=65, in=180] node[text=black, near end, below]{$A$} (0.8,1) to [out=0, in=95] 
            (1,0.5) to [out=270 , in= 20] (0.7,0) to (0.3, 0) to [out = 140, in= 280] (0,0.3);
            \filldraw[Dandelion!85!black] (0,0.3)  (0.8,1) to (1,0.5) to (0.4,0) to  (0,0.3);
            \node[black]at (0.55,0.45) {$J$};

                \node[black, left]at (0,0.3) {$(0,v_1)$};
                    \node[black, below]at (0.4,0) {$(u_1, 0)$};
                        \node[black, above]at (0.8,1) {$(u_2,1)$};
                            \node[black, right]at (1,0.5) {$(1,v_2)$};
        
              \draw[->, xshift=1.6cm] (0,0) -- (1,0) node[below]{$u$} ;
		\draw[->, xshift=1.6cm] (0,0) -- (0,1) node[left]{$v$} ;
        \draw[Black!25, xshift=1.6cm] (1,0) -- (1,1);
        \draw[Black!25, xshift=1.6cm] (0,1) -- (1,1);
        \draw[Blue, thick, xshift=1.6cm, thick] (0, 0.3) -- node[below, midway]{$\eta_1$} (0.4, 0.65);
        \draw[Blue, thick, xshift=1.6cm, thick] (0, 0.3) -- node[below, midway]{$\eta_2$} (0.4, 0);
        \draw[YellowGreen, xshift=1.6cm, thick] (0.4, 0.65) -- node[above, midway]{$\sigma_1$} (0.8, 1);
        \draw[Blue, xshift=1.6cm, thick] (0.8, 1) -- node[right, midway]{$\eta_3$} (1, 0.5);
        \draw[Blue, xshift=1.6cm, thick] (1, 0.5) -- node[above, midway]{$\eta_4$} (0.8, 0.33333);
        \draw[YellowGreen, xshift=1.6cm, thick] (0.8, 0.333333) -- node[above, midway]{$\sigma_2$} (0.4, 0);
        \draw[Sepia, xshift=1.6cm, thick] (0.4, 0.65) node[black, above left]{$(u_1,w_1)$} -- node[above, midway]{$\varsigma$} (0.8, 0.33333) node[black, below right]{$(u_2,w_2)$};
         \draw[xshift=1.6cm, thick] (0.4, 0.65) circle (0.006cm);
         \draw[xshift=1.6cm, thick] (0.8, 0.3333) circle (0.006cm);

    \end{tikzpicture}
    \caption{(Left) An example of set $A$ and subset $J$; we have $\uC_2(A) = \uC_2(J) = 1$. (Right) Support of the copula $C^*$, as per the construction in Appendix~\ref{sect:construction-copula}.}
	\label{fig:quadr}
\end{figure}

   Choose $u_1,u_2,v_1,v_2\in[0,1]$ such that $(u_1,0),(u_2,1),(0,v_1),(1,v_2) \in A$. 
     Let $J$ be the convex hull generated by these four points. From the convexity of $A$, we have $J\subseteq A$; thus, $\uC_2(J) \leq \uC_2(A).$ We will construct a copula $C^*$ supported on the boundary of $J$ and a segment in its interior. We depict an example of $J$ and the support of $C^*$ in Figure~\ref{fig:quadr}. 

We first discuss the case where $u_1,u_2,v_1,v_2\in(0,1)$, meaning none of the four points forming $J$ are corners. By Theorem~\ref{th:rearrangements}, we can reflect each coordinate without changing $\uC_2(J)$, so it is sufficient to only construct $C^*$ for the case where $u_1\leq u_2$.

We divide $[0,1]^2$ into three vertical strips, labelled $K_1,K_2,K_3$, where $K_1$ ranges from $[0,u_1]$, $K_2$ from $(u_1,u_2]$ and $K_3$ from $(u_2,1]$. Let $w_1$ be the height of the upper edge of $J$ at $u_1$, and $w_2$ be the height of its lower edge at $u_2$, as identified in Figure~\ref{fig:quadr}. 
Define
\begin{equation*}
    r_1 = \frac{u_1}{w_1}, \quad r_2 = \frac{1-u_2}{1-w_2}, \quad \text{and}\quad r_3 = 1-r_1 - r_2;
\end{equation*}
If $r_3<0$, interchanging the coordinates and keeping the same construction will yield $r_3^*>0$, as this would mean exchanging $u_1$ and $v_1$ and $u_2$ and $v_2$; it is straightforward to verify this fact algebraically. We thus suppose that $r_3\geq 0$.  Let $\eta_1, \eta_2, \eta_3, \eta_4, \sigma_1,\sigma_2,\varsigma$ be uniform measures on the segments identified in Figure~\ref{fig:quadr}. 

The case $w_1\geq w_2$ is the one illustrated in Figure~\ref{fig:quadr}. For the case $w_1<w_2$, the second-dimension projection of the two segments associated to the measures $\sigma_1$ and $\sigma$ overlap. We thus require $\varsigma_{\rm alt}$ as the uniform distribution on the segment from $(u_1,0)$ to $(u_2,1)$ instead of $\varsigma$ in the mixture. For each case, we construct the copula $C^*$ as a mixture of the measures in the proportions given by Table~\ref{tab:quadr}. 

\begin{table}[tb]
    \centering
    \begin{tabular}{cccc}
         \toprule
         Strip & Measure & Weight in mixture if $w_1\geq w_2$ &  Weight in mixture  if $w_1<w_2$\\
         \midrule
         $K_1$ & $\eta_1$ & $r_1(w_1-v_1)$  &  $r_1(w_1-v_1)$  \\
         & $\eta_2$ & $r_1 v_1$  &  $r_1 v_1$  \\
         $K_2$ & $\sigma_1$ &  $(1-r_2)(1-w_1)$ &   $r_1(1-w_1)$ \\
          & $\sigma_2$ &  $(1-r_1)w_2$ & $r_2w_2$  \\
           & $\varsigma$ &  $r_3(w_1-w_2)$ &   0 \\
         & $\varsigma_{\rm alt}$ &  0 &  $r_3$  \\
         $K_3$ & $\eta_3$ & $r_2(1-v_2)$  &  $r_2(1-v_2)$  \\
          & $\eta_4$ & $r_2(v_2-w_2)$  &  $r_2(v_2-w_2)$  \\
    \bottomrule
    \end{tabular}
    \caption{Weights for the mixture constructing $C^*$. Each measure is a uniform distribution on its associated segment.}
    \label{tab:quadr}
\end{table}

We argue that $C^*$ has uniform marginals.
For the first marginal: The total weight associated to measures on the strip $K_1$ is $r_1(w_1 - v_1) + r_1v_1 = u_1$. Similarly, the total weight of the measures on $K_2$ is $u_2-u_1$ (in both cases) and on $K_3$ is $1-u_2$, which all correspond to the the width of that strip. Each measure is uniformly distributed on its associated segment, and their first-coordinate projection of these segments covers all of their respective strip's width, thus yielding a uniform marginal. 

For the second marginal: For a uniform distribution on a segment, if the projection of the segment has length $h$, then the density of that marginal is $1/h$ over that projection. Let us remark that, for both the cases $w_1\geq w_2$ and $w_1>w_2$, the respective weight of each measure is proportional to the length of the second coordinate projection of its segment. Hence, for $w_2\geq w_1$, the second marginal of $C^*$ has density, for $x\in(0,1)$, 
\begin{equation*}
    f(x) = \underbrace{r_1 \id_{\{x\in(0,w_1)\}}}_{\text{from } \eta_1,\eta_2} + \underbrace{r_2 \id_{\{x\in(w_2,1)\}}}_{\text{from }\eta_3,\eta_4} + \underbrace{(1-r_2) \id_{\{x\in(w_1,1)\}}}_{\text{from }\sigma_1} + \underbrace{(1-r_1) \id_{\{x\in(0,w_2)\}}}_{\text{from }\sigma_2} + \underbrace{r_3 \id_{\{x\in(w_2,w_1)\}}}_{\text{from }\varsigma} = 1.   
\end{equation*}
For the case $w_1<w_2$, instead we have
\begin{equation*}
    f(x) = \underbrace{r_1 \id_{\{x\in(0,w_1)\}}}_{\text{from } \eta_1,\eta_2} + \underbrace{r_2 \id_{\{x\in(w_2,1)\}}}_{\text{from }\eta_3,\eta_4} + \underbrace{r_1 \id_{\{x\in(w_1,1)\}}}_{\text{from }\sigma_1} + \underbrace{r_2 \id_{\{x\in(0,w_2)\}}}_{\text{from }\sigma_2} + \underbrace{r_3\id_{\{x\in(0,1) \}}}_{\text{from }\varsigma_{\rm alt}} = 1.   
\end{equation*}
In both case, this indicates that the second marginal is also a uniform distribution on $(0,1)$. 

Next, for the case where either $u_1 = v_1 = 0$  or $u_2 = v_2 = 1$, meaning one of the  points touched by $A$ is a corner. By applying a reflection using Theorem~\ref{th:rearrangements} if necessary, we can suppose that $u_1 = v_1 = 0$ and $u_2,v_2\in(0,1)$. We construct $C^*$ as a mixture of three uniform distributions on the following segments: 
\begin{align*}
&\text{For segment }(0,0)\text{---}(u_2,1),\text{ assign weight }\frac{u_2(1-v_2)}{1-u_2v_2};\\
&\text{for segment }(0,0)\text{---}(1,v_2),\text{ assign weight }\frac{(1-u_2)v_2}{1-u_2v_2};\\
&\text{for segment }(u_2,1)\text{---}(1,v_2),\text{ assign weight }\frac{(1-u_2)(1-v_2)}{1-u_2v_2}.
\end{align*}
It is straighforward to verify, similarly as above, that the marginals are indeed uniform. 

Finally, for the case where two or more corners are touched, either two of them lie on the same edge, in which case one can use the weights above, or the two corners are $(0,0)$ and $(1,1)$ or $(0,1)$ and $(1,0)$, for which the comonotonic copula and countermonotonic copula are respectively supported on the closed segment from one to the other.  
\end{proof}



\begin{remark}
The segment $\varsigma$ or $\varsigma_{\rm alt}$ is necessary to obtain support of a copula in general. Not every quadrilateral touching the four boundaries of $[0,1]^2$ can support  a copula in its four edges, even by adjusting the weights attributed to (segments of) each of them. To see this more clearly, suppose we have the quadrilateral given by connecting the four points $(\varepsilon,0)$, $(0,1-\varepsilon)$, $(1-\varepsilon, 1)$, and $(1,\varepsilon)$, where $0<\varepsilon<0.25$. Then, the uniform margin requirement implies there is at least $1-2\varepsilon$ units of mass associated to segments $(0,1-\varepsilon)$\textemdash$(1-\varepsilon,1)$ and $(\varepsilon,0) $\textemdash$(1,\varepsilon)$. The second-coordinate projection of these segments has Lebesgue measure $2\varepsilon$; thus, satisfying the uniform margin in the first dimension will inevitably lead to not satisfying it in the other dimension since $1-2\varepsilon>2\varepsilon$ by assumption. 
\end{remark}


{
\small
\begin{thebibliography}{10}

\bibitem[\protect\citeauthoryear{Aharoni et al.}{2002}]{AHKM02}
Aharoni, R., Holzman, R., Krivelevich, M. and Meshulam, R. (2002).
Fractional planks. \emph{Discrete \& Computational Geometry},
\textbf{27}(4), 585--602.

\bibitem[\protect\citeauthoryear{Akopyan et al.}{2019}]{AKP19}
Akopyan, A., Karasev, R. and Petrov, F. (2019). Bang's problem and symplectic invariants. \emph{Journal of Symplectic Geometry}, \textbf{17}(6), 1579--1611. 

\bibitem[\protect\citeauthoryear{Bang}{1951}]{B51}
Bang, T. (1951). A solution of the ``plank problem''. \emph{Proceedings of the American Mathematical Society}, \textbf{2}(6), 990--993.


	\bibitem[\protect\citeauthoryear{Bednarski}{1982}]{B82}
	Bednarski, T. (1982).  Binary experiments, minimax tests and 2-alternating capacities. \emph{Annals of Statistics}, {\bf 10}(1), 226--232.  
	


	\bibitem[\protect\citeauthoryear{Besser and Jorswieck}{2020}]{BJ20}
	{Besser, K. and Jorswieck, E. A.} (2020).  Reliability bounds for dependent fading wireless channels. \emph{IEEE Transactions on Wireless Communications}, {\bf 19}(9), 5833--5845.  
	
	
\bibitem[\protect\citeauthoryear{Blanchet et al.}{Blanchet et al.}{2025}]{BLLW24} Blanchet, J., Lam, H., 
Liu, Y. and Wang, R. (2025). Convolution bounds on quantile aggregation. \textit{Operations Research}, \textbf{73}(5), 2761--2781. 

\bibitem[\protect\citeauthoryear{Boerma et al.}{2025}]{BTZ25}
Boerma, J., Tsyvinski, A. and Zimin, A. P. (2025). Sorting with teams.
\emph{Journal of Political Economy}, \textbf{133}(2), 421--454.



	
\bibitem[\protect\citeauthoryear{Choquet}{Choquet}{1954}]{C54} 
Choquet, G. (1954). Theory of Capacities. \textit{Annales de l'Institut Fourier}, \textbf{5}, 131--295. 

\bibitem[\protect\citeauthoryear{Dempster}{Dempster}{2008}]{D08}  Dempster, A. P. (2008). Upper and lower probabilities induced by a multivalued mapping. In \emph{Classic works of
the Dempster-Shafer theory of belief functions}, 57--72. 

\bibitem[\protect\citeauthoryear{Denuit et al.}{Denuit et al.}{2006}]{DDGKL06}
Denuit, M., Dhaene, J., Goovaerts, M., Kaas, R. and Laeven, R. (2006). Risk measurement with equivalent utility principles. \emph{Statistics and Decisions}, \textbf{24}(1), 1--25.




    

\bibitem[\protect\citeauthoryear{Edwards}{Edwards}{1978}]{E78}
Edwards, D.~A. (1978). On the existence of probability measures with given marginals. \emph{Annales de l'Institut Fourier}, \textbf{28}(4), 53--78.
    
	
\bibitem[\protect\citeauthoryear{Embrechts et al.}{Embrechts et al.}{2013}]{EPR13}  Embrechts, P.,  Puccetti, G. and R\"{u}schendorf, L.  (2013). Model uncertainty and VaR aggregation. {\em Journal of Banking and Finance}, \textbf{37}(8), 2750--2764.

\bibitem[\protect\citeauthoryear{Embrechts et al.}{Embrechts et al.}{2015}]{EWW15} {Embrechts, P., Wang, B. and Wang, R.} (2015). Aggregation-robustness and model uncertainty of regulatory risk measures.  {\em Finance and Stochastics},  \textbf{19}(4), 763--790.



\bibitem[\protect\citeauthoryear{Gaffke and R\"uschendorf}{1984}]{GR84} 
Gaffke, N. and Rüschendorf, L. (1984). On the existence of probability measures with given marginals. \emph{Statistics \& Decisions}, \textbf{2}, 163--174.

\bibitem[\protect\citeauthoryear{Gilboa and Schmeidler}{1989}]{GS89} 
Gilboa, I. and Schmeidler, D. (1989). Maxmin expected utility with non-unique prior. \emph{Journal of Mathematical Economics}, \textbf{18}(2), 141--153.



\bibitem[\protect\citeauthoryear{Hall}{Hall}{1935}]{H35}
Hall, P. (1935). On representatives of subsets. \emph{Journal of the London Mathematical Society}. \textbf{10}(1), 26--30. 

\bibitem[\protect\citeauthoryear{Huber and Ronchetti}{Huber and
  Ronchetti}{2009}]{HR09}
Huber, P.~J. and Ronchetti E.~M.  (2009).
  {\em Robust Statistics.\/} Second ed.,
 Wiley Series in Probability and Statistics. Wiley, New Jersey. 


\bibitem[\protect\citeauthoryear{Huber and Strassen}{Huber and
  Strassen}{1973}]{HS73}
Huber, P.~J. and Strassen V.  (1973). Minimax tests and the Neyman-Pearson lemma for capacities. \emph{Annals of Statistics}, \textbf{1}(2), 252--263.  


\bibitem[\protect\citeauthoryear{Jaffray}{Jaffray}{1991}]{J91}
Jaffray, J.-Y. (1991). Belief functions, convex capacities, and decision making. \emph{Mathematical psychology: Current developments}, 127--134.



\bibitem[\protect\citeauthoryear{Karp}{Karp}{1972}]{K72}
Karp, R.M. (1972). Reducibility among Combinatorial Problems. In: Miller, R.E., Thatcher, J.W., Bohlinger, J.D. (eds) \textit{Complexity of Computer Computations.} The IBM Research Symposia Series. Springer, Boston. 

\bibitem[\protect\citeauthoryear{Keevash and Lim}{2026}]{KL26}
Keevash, P. and Lim, J. (2026). On the largest sum-free subset of the lattice cube.
arXiv preprint arXiv:2605.00816.

\bibitem[\protect\citeauthoryear{Kellerer}{Kellerer}{1984}]{K84}
Kellerer, H. G. (1984). Duality theorems for marginal problems. \emph{Zeitschrift für Wahrscheinlichkeitstheorie und verwandte Gebiete}, \textbf{67}, 399--432.

\bibitem[\protect\citeauthoryear{Lepsveridze and Sun}{2026}]{LS26}
Lepsveridze, S. and Sun, Y. (2026). Size of the largest sum-free subset of $[n]^3$ and $[n]^4$.
\emph{International Mathematics Research Notices}, \textbf{2026}(9), rnag081.



\bibitem[\protect\citeauthoryear{Marinacci and Montrucchio}{Marinacci and Montrucchio}{2004}]{MM04}
{Marinacci, M. and Montrucchio L.} (2004). Introduction to the mathematics of ambiguity.  In \emph{Uncertainty in Economic Theory}, I. Gilboa, Ed.  Routledge, New York, 46--107.




\bibitem[\protect\citeauthoryear{McNeil et al.}{McNeil et al.}{2015}]{MFE15}
{McNeil, A. J., Frey, R. and Embrechts, P.} (2015). \emph{Quantitative
Risk Management: Concepts, Techniques and Tools}. Revised Edition.  Princeton, NJ:
Princeton University Press.




\bibitem[\protect\citeauthoryear{Puccetti et al.}{2019}]{PRWW19}
{Puccetti, G., Rigo, P., Wang, B. and Wang, R.} (2019).
	Centers of probability measures without the mean.
 {\em Journal of Theoretical Probability}, \textbf{32}(3), 1482--1501.

\bibitem[\protect\citeauthoryear{Puccetti and Wang}{Puccetti and Wang}{2015}]{PW15}
Puccetti, G. and Wang R.  (2015).
Extremal dependence concepts. \emph{Statistical Science},  \textbf{30}(4), 485--517.


\bibitem[\protect\citeauthoryear{Ramachandran and R\"uschendorf}{1995}]{RR95}
Ramachandran, D. and R\"uschendorf, L.  (1995). A general duality theorem for marginal problems. \emph{Probability Theory and Related Fields}, \textbf{101}, 311--319.
 


\bibitem[\protect\citeauthoryear{Roberson}{2006}]{R06}
Roberson, B.  (2006). The Colonel Blotto game. \emph{Economic Theory}, \textbf{29}(1), 1--24.
 
 
 \bibitem[\protect\citeauthoryear{Santambrogio}{Santambrogio}{2015}]{S15} Santambrogio, F. (2015). \emph{Optimal Transport for Applied Mathematicians}. Springer, New York.
 
\bibitem[\protect\citeauthoryear{Schmeidler}{Schmeidler}{1989}]{S89}
{Schmeidler, D.} (1989). Subjective probability and expected utility without additivity.
\emph{Econometrica}, \textbf{57}(3), 571--587.

 \bibitem[\protect\citeauthoryear{Shaked and Shanthikumar}{Shaked and Shanthikumar}{2007}]{SS07} Shaked, M. and Shanthikumar, J. G. (2007). \emph{Stochastic
Orders}. Springer Series in Statistics.

 \bibitem[\protect\citeauthoryear{Shortt}{Shortt}{1986}]{S86}
Shortt, R. M. (1986). The singularity of extremal measures. \emph{Real Analysis Exchange}, \textbf{12}, 205--215.

 
\bibitem[\protect\citeauthoryear{Siegel}{2009}]{S09}
Siegel, R. (2009). All‐pay contests. \emph{Econometrica}, \textbf{77}(1), 71--92.

\bibitem[\protect\citeauthoryear{Skala}{1993}]{S93}
Skala, H. J. (1993). The existence of probability measures with given marginals. \emph{The Annals of Probability}, \textbf{21}(1), 136--142.

\bibitem[\protect\citeauthoryear{Sklar}{1959}]{S59}
Sklar, M. (1959) Fonctions de répartition à N dimensions et leurs marges. \emph{Annales de l’Institut de Statistiques de l'Université de Paris}, \textbf{8}(3), 229--231.


\bibitem[\protect\citeauthoryear{Strassen}{1965}]{S65}
Strassen, V. (1965). The existence of probability measures with given marginals. \emph{The Annals of Mathematical Statistics}, \textbf{36}(2), 423--439.
 
 
\bibitem[\protect\citeauthoryear{Vovk et al.}{Vovk et al.}{2022}]{VWW22}
Vovk, V., Wang, B. and Wang, R. (2022). Admissible ways of merging p-values
  under arbitrary dependence.  \emph{Annals of Statistics}, \textbf{50}(1), 351--375. 

\bibitem[\protect\citeauthoryear{Vovk and Wang}{Vovk and Wang}{2020}]{VW20}
Vovk, V. and  Wang, R. (2020).
Combining p-values via averaging.  \emph{Biometrika},  \textbf{107}(4), 791--808.
 

 \bibitem[\protect\citeauthoryear{Wang and Wang}{Wang and Wang}{2016}]{WW16}
Wang, B. and  Wang, R. (2016).
 Joint mixability.  \emph{Mathematics of Operations Research}, \textbf{41}(3), 808--826.

\bibitem[\protect\citeauthoryear{Wang, Peng and Yang}{Wang
  et~al.}{2013}]{WPY13}
Wang, R., Peng, L. and Yang, J. (2013).
 Bounds for the sum of dependent risks and worst Value-at-Risk with
  monotone marginal densities.
 {\em Finance and Stochastics}, \textbf{17}(2), 395--417.
 
 
\bibitem[\protect\citeauthoryear{Wasserman and Kadane}{Wasserman and Kadane}{1990}]{WK90}
  Wasserman, L. A. and Kadane, J. B. (1990). Bayes' theorem for Choquet capacities. \emph{Annals of Statistics}, \textbf{18}(3), 1328--1339.

\bibitem[\protect\citeauthoryear{Wolfenson and Fine}{Wolfenson and Fine}{1982}]{WF82}
  Wolfenson, M. and Fine, T. L. (1982). Bayes-like decision making with upper and lower probabilities. \emph{Journal of the American Statistical Association}, \textbf{77}(377), 80--88.

\end{thebibliography}
}
\end{document}